\documentclass[12pt]{article}
\usepackage[margin=1in]{geometry}
\usepackage{mathtools}
\usepackage{amsmath}
\usepackage{amsthm,bm}
\usepackage{textcomp}
\usepackage{graphicx}
\usepackage{amssymb}
\mathtoolsset{showonlyrefs=true}
\usepackage{float}
\usepackage{pdflscape}
\usepackage{pgfplots}
\usepackage{mathrsfs}
\usepackage[numbers]{natbib}    
\usepackage{hyperref}
\usepackage{bbm}
\usepackage{comment}
\numberwithin{equation}{section}
\usepackage{xcolor}
\usepackage[normalem]{ulem}
\DeclareMathOperator*{\argmax}{argmax}
\newtheorem{thm}{Theorem}[section]

\newtheorem{prop}[thm]{Proposition}

\theoremstyle{definition}

\usepackage{bbm}

\theoremstyle{remark}
\newtheorem{rem}[thm]{Remark}

\usepackage{graphicx, bigints}
\usepackage{graphicx}

\usepackage{graphicx}    
\usepackage{caption}    
\usepackage{subcaption}  

\NewDocumentCommand{\citeayn}{m}{%
  \citeauthor{#1}\ (\citeyear{#1}),\ \citenum{#1}%
}
\RenewDocumentCommand{\citep}{o o m}{%
  [%
    \citeayn{#3}%
    \IfValueT{#1}{, #1}%
    \IfValueT{#2}{ #2}%
  ]%
}

\pgfplotsset{compat=1.18}
\begin{document}
\title{Fare Game: A Mean Field Model of Stochastic \\Intensity Control in Dynamic Ticket Pricing
}

\makeatletter
  \let\origthanks\thanks
  \renewcommand{\thanks}[1]{
    \origthanks{
      \fontsize{9pt}{15pt}\selectfont#1}
  }
\makeatother

\author{Burak Ayd{\i}n\thanks{Department of Operations Research \& Financial Engineering, Princeton University, {\em ba7735@princeton.edu}}\and Emre Parmaks{\i}z\thanks{Department of Mathematics, Massachusetts Institute of Technology, {\em yep@mit.edu}}\and Ronnie Sircar\thanks{Department of Operations Research \& Financial Engineering, Princeton University, {\em sircar@princeton.edu}}}
\date{December 30, 2025}
\maketitle
\begin{abstract}
    We study the dynamic pricing of discrete goods over a finite selling horizon. 
    One way to capture both the elastic and stochastic reaction of purchases to price is through a model where sellers control the intensity of a counting process, representing the number of sales thus far. The intensity describes the probabilistic likelihood of a sale, and is a decreasing function of the price a seller sets. A classical model for ticket pricing, which assumes a single seller and finite time horizon, is by Gallego and van Ryzin (1994) and it has been widely utilized by airlines, for instance. Extending to more realistic settings where there are multiple sellers, with finite inventories, in competition over a finite time horizon is more complicated both mathematically and computationally. We introduce a dynamic mean field game of this type, and some numerical and existence results. In particular, we analyze the associated coupled system of Hamilton-Jacobi-Bellman and Kolmogorov differential-difference equations, and we prove the existence and uniqueness results under certain conditions.  Then, we demonstrate a numerical algorithm to find this solution and provide some insights into the macroeconomic market parameters. Finally, we present a qualitative comparison of our findings with airfare data.
\end{abstract}

\section{Introduction}

Deregulation of the airline industry in 1978 ignited the field of revenue management, which was initially known as yield management. Faced with new competitive dynamics, airlines adopted inventory control and dynamic pricing to optimize revenue. These techniques have since diffused across numerous industries, including transportation, entertainment, hospitality, and retail. This paper proposes dynamic pricing strategies for discrete goods such as tickets for an event, in a competitive market over a finite time horizon corresponding to the date of the event.

Sales are modeled as Poisson-type random arrivals of customers to a representative seller.
The probability (or intensity) of these arrivals is decreasing
in the individual seller's price, and it is increasing with the average market price at any given time.
Finite player nonzero-sum stochastic differential games of this type are analytically challenging, even with just two firms, as shown in \citep{Ledvina2011}. They become computationally intractable with more than five firms. To obtain a computationally feasible approximation that yields interpretable results, we propose a mean field game model, which assumes a continuum of firms across the inventory space.

A pioneering work by \citep{Belobaba1987} underscored how initial airline seat inventory practices relied heavily on human judgment to balance discounted and full-fare tickets for revenue maximization. Recognizing the need for systematic approaches, the Expected Marginal Seat Revenue model was subsequently introduced \citep{Belobaba1989}, which significantly advanced revenue management by incorporating probabilistic decision-making for optimal booking limits. These analytical methodologies demonstrated substantial economic benefits. With advances in computational power, dynamic pricing grew prominent across industries, as extensively reviewed by \citep{McGillVanRyzin1999} in the areas of forecasting, overbooking, seat inventory control, and pricing. 

Dynamic pricing and inventory management under stochastic demand were rigorously formalized in \citep{Gallego1994}, where they consider a firm operating in a market with imperfect competition, e.g. a monopoly, such that the firm has the ability to influence its demand by changing the price of a product. The demand in their model is given by a Poisson-type process with intensity $\lambda(p)$, which is non-increasing in price $p$. At every arrival, one unit is sold at the quoted price, and inventory is decreased by one, until it reaches zero. Their framework illustrated the efficacy of dynamic pricing policies.

A related direction is discrete-choice models introduced by \citep{TalluriVanRyzin2004}, which explicitly incorporated consumer substitution behavior, allowing more realistic and profitable pricing strategies. Inventory considerations within dynamic pricing contexts were extensively reviewed by \citep{ElmaghrabyKeskinocak2003}, providing critical insights applicable across airlines, hospitality, and retail sectors. This broadened the applicability of dynamic pricing methodologies beyond their traditional scopes. A more recent survey on dynamic pricing (and learning) can be found in \citep{denBoer}.

Recent research on oligopolistic airline markets, such as \citep{Hortacsu2022}, extended earlier monopolistic frameworks to duopolies and multi-firm setups, highlighting strategic interactions under dynamic pricing with finite inventory and selling horizons. They manage to show improvements in output at the expense of lower welfare levels due to dynamic pricing. Integrating machine learning and game theory into revenue management has also become a prominent area of development. Various methods of game-theory-based reinforcement learning in oligopolistic contexts are analyzed in \citep{CollinsThomas2012}, and comparisons between reinforcement learning and data-driven programming are discussed in \citep{Lange2025}. However, as shown in \citep{Bondoux}, although these reinforcement learning models may be proven to converge to the optimal solution in the market, in practice they require very large amounts of data. 
Due to limitations with scarce data, we instead introduce a model-based game theoretic approach.

Another branch of literature that is of interest to our paper is stochastic intensity control, which has also gained considerable attention parallel to revenue management research. Recognizing the need for stochastic intensity control models in perishable goods pricing, earlier models were generalized to allow nonhomogeneous but deterministic Poisson consumer arrival rates in \citep{ZhaoZheng2000}, significantly improving revenue outcomes compared to static approaches. Competitive revenue management under deterministic demand arrival rates was analyzed by \citep{GallegoHu2014}, who studied dynamic pricing equilibrium through shadow prices and extended results to stochastic intensity-controlled games. Similarly, intensity-controlled stochastic exploration in dynamic Cournot games was studied in \citep{LudkovskiSircar}, revealing strategic adjustments firms make under uncertainty to manage exhaustible resources. These insights are applicable to stochastic demand pricing scenarios. 

Despite extensive studies on revenue management, intensity control, and their intersections, mean-field games of intensity control for markets characterized by nonhomogeneous consumer arrival rates remain unexplored. Our paper fills this gap, providing a novel analysis of competition through mean-field game theory as an alternative to analytically and numerically challenging finite-horizon multi-firm competitive models. Mean-field game (MFG) theory is an active research area, thoroughly covered in texts such as \citep{CarmonaDelarue2018}. It has been applied to various problems related to ours, including Bertrand and Cournot mean-field games \citep{ChanSircar2015}, optimal exploitation of exhaustible resources \citep{GraeweHorstSircar2022}, oil markets and fracking \citep{ChanSircar2017}, optimal execution in a dynamic pricing context \citep{Donnelly}, portfolio liquidation problems \citep{FuGraeweHorstPopier2021}, or cryptocurrency mining \citep{LiReppenSircar2022}. In the context of the latter model, we refer also to the recent existence results for MFGs of intensity control in \citep{Garcia_crypto}. 
While our model is related to this paper, where the authors develop a general existence theory for mean field games of intensity control in both discrete and continuous time, some of the differences can be summarized as follows. Their analysis is for games with continuous state spaces, it uses randomized (mixed) strategies and a Kakutani fixed point theorem, while we study a setting where the state space is discrete, with pure strategies, and use a Schauder fixed point theorem.

The remainder of the paper is organized as follows. Section \ref{sec:2} introduces the detailed model and its associated coupled Hamilton-Jacobi-Bellman and Kolmogorov differential-difference equations that characterize the mean field market equilibrium. Our main existence and uniqueness theorems are given in Section \ref{sec:3}, with their proofs in Appendix \ref{app:4}. Section \ref{sec:numerics} outlines the numerical algorithm, shows our results, discusses the algorithm’s numerical stability and effectiveness, and presents a qualitative comparison of our model with airfare data. Finally, Section \ref{sec:concs} concludes 
and discusses directions for future research.


\section{Stochastic Intensity Mean Field Game} \label{sec:2}

We consider a continuum of firms distributed across our state space, comprising positive integers. The states represent the number of remaining items, such as tickets, in each firm's inventory. The model encompasses two primary interactions. First, we have the individual optimization of the firms in response to average market prices, strategically adjusting their prices based on their competitors' actions. Second, we have the dynamics of the market itself, which are shaped by the pricing strategies adopted by each firm. The general market dynamics and the actions of individual terms are interdependent through the customer arrival intensity function. 

A mean field equilibrium is a set of strategies that obey the following: assuming a certain measure flow (evolution) of the system, every representative agent, given their state, chooses their optimal adaptive control (strategy). The resulting controls cause the system to evolve according to Kolmogorov’s equation, which in turn produces a new measure flow. If the input and output measure flows coincide, we are in mean field equilibrium (\citep{CarmonaDelarue2018}, \citep{Lacker2018MFGNotes}). To put it in another way, a mean field equilibrium is reached when the average of individual firm best responses to a posited mean price, coincides with that mean price (at every time). In other words there is an optimization coupled with a fixed point. 

\subsection{Model Dynamics} 
We denote the distribution of firms' inventories by $m(t) = \{m_k(t)\}_{k=0}^{K}$, where $m_k(t)$ is the proportion of the firms with $k$ items in their inventory at time $t$, and $K < \infty$ denotes the maximum possible inventory level. We have $\sum_{k=0}^{K} m_k(t)= 1$ for all times $0\leq t \leq T$, where $T<\infty$ is the event horizon (for instance time of flight or concert). Firms have finite inventories and are incentivized to sell, as any unsold inventory holds no residual value at the end of the selling period. The initial inventory distribution is given by $m_k(0)=M_k$.

If we denote a firm's inventory by $X$, its dynamics are 
\begin{equation}\label{eq:2.1}
dX_t = - \mathbbm{1}_{\{X_t > 0\}}dN_t^{\lambda},
\end{equation}
where $N_t^{\lambda}$ is a counting process with stochastic and time-varying intensity $\lambda_t \geq 0$, meaning
\begin{align}\label{eq:lambda_prob}
    \mathbbm{P}(N^\lambda_{t+h}-N^\lambda_t=1)= \lambda_th+o(h),\qquad
    \mathbbm{P}(N^\lambda_{t+h}-N^\lambda_t\ge2)=o(h).
\end{align}
See for instance \citep{bremaud1981}.
The intensity $\lambda_t$ is a function that depends on the price $p_t$ quoted by the firm, the average price $\bar{p}(t)$ in the market, and the firm inventory distribution $m$ on the state space:
\begin{equation}\label{eq:2.2}
\lambda_t = \lambda (p_t, \bar{p}(t), m(t)).
\end{equation}
The function $\lambda$ is assumed decreasing in the firm's own price $p_t$, so decreasing price increases the probability of a sale.
It is increasing in the average price in $\bar{p}(t)$, so with a higher average price, a firm will expect more demand. 

A representative seller optimizes expected discounted revenue and the value functions $V(t)= \{ V_k(t) \}_{k=0}^{K}$ when they have $k$ tickets at time $t\leq T$
for each $k$ are given by
\begin{equation}\label{eq:2.3}
V_{k}(t) = \sup_{{p_s \geq 0}} \mathbbm{E} \left\{  \int_{t}^{T} e^{-r(s-t)} p_s \mathbbm{1}_{\{X_s > 0\} }\,dN_s^\lambda \mid X_t = k\right\},
\end{equation}
where the supremum is taken over adapted 
price trajectories $p_s$ over the trading period $t \leq s \leq T$. 

For $k\in \mathcal{K} := \{1, \dots , K\}$, the associated HJB differential-difference equation is 
\begin{equation}\label{eq:2.4}
    \dfrac{dV_k}{dt} - rV_k + \sup_{p\geq 0} \Big\{ \lambda(p,\bar{p}, m) (p - \Delta V_k)\Big\} = 0,
\end{equation}
where we define the discrete difference $\Delta V_k := V_k - V_{k-1}$. We have the terminal and boundary conditions
\begin{equation}\label{eq:2.5}
        V_k(T) = 0, \quad k\in\mathcal{K},\qquad
        V_0(t) = 0,\quad  \ t \in[0,T],
\end{equation}
where the first condition means that the remaining goods in the inventory at the end of the selling period will remain unsold and thus generate no value. The second condition means for a firm starting the trading period with zero inventory, the value function is identically zero.

For simplicity, we will consider smooth $\lambda(p, \bar{p}, m)$, and assume that there is a unique supremum in \eqref{eq:2.4}, 
which is expressed as
\begin{equation}\label{eq:2.6}
    p_k^* = \argmax_{p\geq0} \Big \{ \lambda(p,\bar{p}, m) (p - \Delta V_k) \Big \}.
\end{equation}
We define the mean-field equilibrium average price in the market at time $t$ as:
\begin{equation}\label{eq:2.9}
    \bar{p}(t) = \dfrac{1}{\eta(t;m)} \sum_{k = 1}^{K} m_k(t) p_k^*(t),
\end{equation}
where 
\begin{equation}\label{etadef}
    \eta(t;m) := \sum_{k = 1}^{K} m_k(t) = 1 - m_0(t)
\end{equation}
 is the proportion of active firms. It is important to note that the inactive sold-out firms are omitted from the average calculation, as $\{k=0\}$ is an absorbing state. We need to keep track of the inactive firms as they do not affect the average price.
 
The second main equation for our system, known as the Kolmogorov equation, describes the time-evolution of the distribution $m(t)$ when the firms set prices $\{p_k^{*}(t)\}_{k=1}^{K}$:
\begin{equation}\label{eq:2.10}
\begin{split}
& \dfrac{dm_0}{dt} = \lambda(p_{1}^*, \bar{p}, m) m_{1}(t), \\
& \dfrac{dm_k}{dt} = \lambda(p_{k+1}^*, \bar{p}, m) m_{k+1}(t) - \lambda(p_{k}^*, \bar{p}, m) m_{k}(t), \qquad k \in \mathcal{K}_{-1} \\
& \dfrac{dm_{K}}{dt} = -\lambda(p_{K}^*, \bar{p}, m) m_{K}(t).
\end{split}
\end{equation}
subject to an initial condition $m_k(0) = M_k$, and $\mathcal{K}_{-1} :=\{1,\dots, K-1\}$. For example, the second line of \eqref{eq:2.10} describes how $m_k$ increases if a firm with inventory $k+1$ receives a demand, therefore its inventory level reduces to $k$, or how it decreases if a $k$-inventory firm receives a demand and drops to $k-1$ inventory.

In the literature, for Brownian setup, there are adaptations of the standard optimal control verification results \citep[Theorem 9.12]{Lacker2018MFGNotes}, \citep[Proposition 7.5]{CarmonaDelarue2018}. These can be mimicked in our case using results in \citep{bremaud1981}. In our setting, it is  necessary to use Itô's lemma and Feynman--Kac formulas adapted to jump processes to mimic the standard proof of the verification theorem (see, for instance, \citep{Zhu2015}, \citep{AndreaPDE}).

The HJB and Kolmogorov equations together characterize the mean field equilibrium in the market. Their coupling through the intensity function $\lambda(p, \bar{p}, m)$ allows us to analyze the interplay between firm strategies, market dynamics, and the resulting equilibrium as the system evolves over time. 
\subsection{Linear Intensity Model} \label{sec:3.0}
In the remainder of the paper, we work with the piecewise linear intensity function
\begin{equation}\label{eq:3.0}
    \lambda(p, \overline{p},m) = (a(m) - p + c(m)\overline{p})^{+},
\end{equation} 
where 
\begin{equation}\label{eq:3.0b}
        a(m) = \dfrac{1}{1 + \epsilon \sum_{k = 1}^{K} m_k}=\dfrac{1}{1+\epsilon\eta(m)}, \qquad
        c(m) = \dfrac{\epsilon\eta(m)}{1+\epsilon\eta(m)},
\end{equation}
and we have suppressed the $t$-dependence of $m$ and $\eta$ in the notation. 
We observe that $a(m) = 1 - c(m)$. The parameter $\epsilon \geq 0$ represents the level of competition in the market. When $\epsilon = 0$, $\lambda(p, \overline{p},m)$ becomes independent of $m$ and $\overline{p}$ and the firms behave like monopolies.

The function in \eqref{eq:3.0} satisfies $\partial_{p}\lambda \leq 0$ and $\partial_{\overline{p}}\lambda \geq 0$. Moreover, when all firms are inactive ($m_0 = 1,\ m_k=0 \ \forall k \in \mathcal{K}$), we have $\eta(m)=0$, and we observe that $\lambda(p, \overline{p},m) = (1-p)^+$ depends only on $p$. This shows how a firm becomes closer to a monopoly as the other firms exit the market by exhausting their inventories.
We refer to \citep{ChanSircar2015} and  \citep[Chapter 6]{vives} for motivation to consider the piecewise linear intensity function. There they consider finitely many agents that have a linear price-demand relation, and by solving the system of equations relating to the price-demand equilibrium, they obtain a relation of the form \eqref{eq:3.0} in the continuum limit $N \to \infty$. The coefficients in \eqref{eq:3.0b} are derived analogously to the classical specification by economists in oligopoly theory, as described for instance in \citep[Chapter 6]{vives} for Cournot and Bertrand games with substitutable goods. In these and our case, $a$ and $c$ represent the maximum demand (at zero prices), and the extent  of competition between firms, respectively. These change with the distribution $m$ through the fraction of sellers who actually receive demand (and so influence the market, as opposed to being effectively redundant). A brief discussion is given in Appendix \ref{app:3}.

Using the expression for $\lambda$ in \eqref{eq:3.0}, and suppressing $a \text{ and } c$'s dependence on $m$, we obtain
\begin{equation} \label{eq:3.3}
    p_k^* = \dfrac{1}{2} \Big ( a + c\overline{p} + \Delta V_k \Big )^{+}.
\end{equation}
As demonstrated in the proof our main theorems, given in Appendix \ref{app:4}, the quantity between the parentheses in \eqref{eq:3.3} is non-negative in equilibrium for small enough $\epsilon$, so we can remove the positive part. Then we obtain a simplified expression for the average price $\bar{p}$:
\begin{equation} \label{eq:3.4}
    \overline{p} = \dfrac{1}{2 - c} \Big ( a + \dfrac{1}{\eta} \sum_{k=1}^{K} m_k \Delta V_k \Big ).
\end{equation}
Plugging these expressions for $p^*_k$ and $\bar{p}$ into \eqref{eq:2.4} and \eqref{eq:2.10}, and defining the function \begin{equation} \label{eq:3.5}
    \phi_\epsilon(V, m) := \dfrac{1}{2 + \epsilon \eta(m)} \Big ( \sum_{k=1}^{K} m_k \Delta V_k - \eta(m) \Big ),
\end{equation}
we arrive at the coupled system of differential equations that together govern the dynamics of the model:
\begin{align} \label{eq:3.6}
\begin{split}
    & \dfrac{dV_k}{dt} - rV_k + \dfrac{1}{4} \Big ( 1 - \Delta V_k + \epsilon \phi_\epsilon(V, m)\Big )^2 = 0,
    \qquad k \in \mathcal{K}  \\
    & \dfrac{dm_k}{dt} - \left[\dfrac{1}{2} \Big(1 - \Delta V_{k+1} + \epsilon \phi_\epsilon(V, m) \Big ) m_{k+1}
    -\dfrac{1}{2} \Big (1 - \Delta V_k + \epsilon \phi_\epsilon(V, m)  \Big )m_k \right] = 0, 
    \quad k \in \mathcal{K}_{-1} \\
     & \dfrac{dm_0}{dt} - \dfrac{1}{2} \Big ( 1 - \Delta V_1 + \epsilon \phi_\epsilon(V, m) \Big ) m_1  = 0, \\
    & \dfrac{dm_K}{dt} + \dfrac{1}{2} \Big ( 1 - \Delta V_K + \epsilon \phi_\epsilon(V, m) \Big ) m_K  = 0, \\
    & V_k(T) = 0, \  \quad k \in \mathcal{K},  \qquad
m_k(0) = M_k,\  \quad k \in \mathcal{K}. 
\end{split}
\end{align} 
A direct observation shows that when the interaction parameter $\epsilon$ is identically zero, the two main equations relating $V$ and $m$ are decoupled. We argue that in the case of non-zero but sufficiently small competition $(\epsilon > 0)$, we can prove existence and uniqueness by linearizing the system around the monopoly case and using the Schauder Fixed Point Theorem \citep[Corollary~11.2, p.~280]{gilbarg2001elliptic}. Our proofs are related to those for the continuous-space Cournot exhaustible resources oil extraction MFG model in \citep{GRABER2023816}.

\subsection{Main Results} \label{sec:3}
The following proposition will be useful in the proofs of the main theorems. 
\begin{prop}\label{prop:2}
    For $\epsilon = 0$, a solution $(V^{(0)}, m^{(0)})  = (\{ V_k^{(0)} \}_{k=1}^{K} , \{ m_k^{(0)} \}_{k=1}^{K})$ to the system of differential equations \eqref{eq:3.6} exists and is unique. Further, for all $k \in \mathcal{K}$, and $t \in [0,T]$, we have
\begin{align}\label{zeroprops}
\begin{split}
    & 0 \leq \Delta V^{(0)}_k(t) \leq 1 -2\rho < 1 \\
    & \dfrac{d\Delta V^{(0)}_k}{dt}(t) \leq 0 \\
    & 0 \leq \eta^{(0)}(t;m^{(0)}) \leq 1 \\
    & 0 \leq m^{(0)}_k(t) \leq 1,
\end{split}
\end{align}
where $\rho = \sqrt{r^2 + r} - r \in [0,\frac12)$. Consequently, for $\lambda^{(0)}_k(t) = \dfrac{1}{2} (1 - \Delta V^{(0)}_k
(t))$, we have
\begin{equation*}
    \rho \leq \lambda^{(0)}_k (t) \leq \dfrac{1}{2}.
\end{equation*}
\end{prop}
The proof is given in Appendix \ref{app:2}.
The bounds \eqref{zeroprops} are in line with our expectations for a solution that is economically sensible. The parameter $\lambda_k^{(0)}$ represents customer arrival rates, hence it should be non-negative as depicted. Similarly, the variables $m_k$ and $\eta$ signify proportions, and thus it is crucial for these values to lie within the range $[0, 1]$ to align with economic rationale. A noteworthy detail is that the $\lambda_k^{(0)}$ has a lower bound of $\rho > 0$, uniformly across $k$ and in $t$.

Our two main theorems are the following:

\begin{thm} \label{thm:2.2}
    In the setting described above, for any initial inventory distribution $\{ M_k\}_{0\leq k\leq K}$, interest rate $r > 0$ and time horizon $T > 0$, there exists a solution $\{V_k(t), m_k(t)\}_{0\leq k \leq K}$ to the HJB-Kolmogorov system of equations if the interaction parameter $\epsilon > 0$ is less than some constant $C_1$, which depends on the model parameters $T, \ K, \ r$. 
\end{thm}

\begin{thm}\label{thm:2.3}
    In the setting described above, for any initial inventory distribution $\{ M_k\}_{0\leq k\leq K}$, interest rate $r > 0$ and time horizon $T > 0$, there exists an upper bound $C_2$ such that if $\epsilon < C_2$, then the solution exists and is unique.
\end{thm}

\begin{rem} \label{rem:2.4}
        The constants appearing in Theorems \ref{thm:2.2} and \ref{thm:2.3} are given explicitly by
    \begin{equation} \label{eq:c1c2}
        C_1(K,T,r) = \frac{1}{K^2 T \, 2^{K+8}} \min \Big\{ 1 + KT, \ \sqrt{r} \Big\},
        \qquad 
        C_2(K) = \frac{2}{K-1}.
    \end{equation}
    Their derivations are provided in the Appendix.  
\end{rem}

The proofs of these theorems are in Appendix \ref{app:4}. For the proof of existence, we consider first-order deviations of the solution $(V,m)$ to \eqref{eq:3.6}, from the monopolist solution $(V^{(0)}, m^{(0)})$, prescribe a compact set as their domain $M$, and derive the differential equations that control them. We derive the integral operators whose fixed points give the solutions to the differential equations.

To achieve these fixed points, one typically applies either the Banach fixed point theorem, to show contraction (as in \citep{GomesMohrSouza2012}), or the Schauder fixed point theorem, which relies on convexity and compactness to guarantee existence. In the literature (see \citep{Lacker2018MFGNotes} and others), both approaches appear, with the Banach route usually requiring the time horizon to be sufficiently small. In our concrete problem, $T$ is fixed by the date of the event. We instead rely on Schauder’s fixed point theorem for sufficiently small values of the interaction parameter $\epsilon$.

In general formulations, assuming certain Lipschitz and regularity conditions on the coefficients, Schauder’s theorem can be applied directly. However, at that level of generality, the assumptions impose restrictions that do not hold in our setting. Typical requirements include Brownian diffusion, joint continuity of the interaction parameter in both the control (price) and measure (state distribution under the Wasserstein-1 metric), and uniformly bounded second derivatives in both variables. Because our model does not satisfy these regularity assumptions, we adapt the proof strategy to our system (Appendix \ref{app:4}).

For uniqueness, we utilize an energy identity. Assuming two solutions to \eqref{eq:3.6}, we prove under restrictions on $\epsilon$ that the solutions are indeed the same. These theorems highlight the existence and structure of equilibrium solutions, laying a solid foundation for subsequent numerical investigations in Section \ref{sec:numerics}. 
\begin{rem}
     In proving the existence of an equilibrium, we aimed to demonstrate the proof strategy rather than achieving the best possible bound for the interaction parameter for our specific intensity function. To that end, for the sake of brevity in calculations, we used simple supremum-type bounds as opposed to more advanced Hölder-type inequalities in our nested integrations, and as such, the bound on the interaction parameter is not optimal, and can be improved from the values in \eqref{eq:c1c2}. We note that it is quite common that mean field game existence and uniqueness results are proven under some small parameter restriction, for instance small level of competition in \citep{GRABER2023816}, or small time horizon in \citep[Chapter 8]{Lacker2018MFGNotes}. However, our numerical simulations indicate that the interaction parameter need not be within that range to achieve existence and uniqueness. We defer obtaining a tighter bound on the interaction parameter as well as generalizing to other intensity functions for future work.
\end{rem}
\begin{rem}
    Comparing the discrete-space proofs to their continuous-space analogs in the literature, the uniqueness proof would follow the same structure in a continuous-space mean field game, relying similarly on an energy identity and monotonicity arguments. The main difference arises in the existence proof: while our discrete-state setting allows us to use compactness arguments in finite dimensions, in the continuous-space case one would need to work with the Wasserstein distance on probability measures and establish additional regularity results.
\end{rem}

\section{Numerical Results}\label{sec:numerics}
In this section, our primary goal is to establish a numerical method for computing the mean-field equilibrium. Subsequently, we will demonstrate the convergence of our theoretical findings to the numerical results in the limit $\epsilon \to 0$. Finally, we will conclude this section by conducting an in-depth analysis of the numerical stability of our solution strategy, explore various ways to establish theoretical convergence proofs, and present qualitative comparison with airfare data.

\subsection{Algorithm}
As discussed in Section \ref{sec:2}, the mean-field equilibrium characterizes a state in which the aggregate behavior of the continuous set of players reaches a fixed point of an optimization problem. In this equilibrium, each player autonomously adopts a Markovian strategy to optimize their individual value function while taking into account the collective behavior of the entire population. As a result, these individually optimized strategies jointly contribute to the overall movement of the population as initially assumed. Our algorithm for finding this fixed point is as follows: 

\subsubsection*{Step 0: Initial guess}
Begin with an initial guess $(\bar{p}^{(0)}(t), m^{(0)}(t))$, and for $n = 0, 1, 2, ...$ do the following steps:

\subsubsection*{Step 1: Solving the HJB equation}
We discretize the HJB equation in time. Using the initial guess $(\bar{p}^{(n)}(t), m^{(n)}(t))$, we employ the Runge-Kutta method for ODEs of order 5(4) to solve.
\begin{equation}
    \dfrac{dV_k}{dt} - rV_k + \sup_{p} \Big\{ \lambda(p,\bar{p}^{(n)}(t), m^{(n)}(t)) (p - \Delta V_k)\Big\}= 0, \quad k\in\mathcal{K} \nonumber
\end{equation}
with the terminal condition $V_k (T) = 0$.
Using the solution, find the prices $p_{k}^{*(n)}(t)$ and intensities $\lambda_{k}^{*(n)}(t) := \lambda(p_{k}^{*(n)}(t), \bar{p}^{(n)}(t), m^{(n)}(t))$. 
\subsubsection*{Step 2: Solving the Kolmogorov equation}
Using intensities in step 1, find $m^{(n+1)}(t)$ by solving the Kolmogorov equation:
\begin{equation}
    \begin{split}
        \dfrac{dm_{0}^{(n+1)}}{dt} &= \lambda_{1}^{*(n)}(t) m_{1}^{(n+1)} \\
        \dfrac{dm_{k}^{(n+1)}}{dt} &= -\lambda_{k}^{*(n)}(t) m_{k}^{(n+1)} + \lambda_{k+1}^{*(n)}(t) m_{k+1}^{(n+1)}, \quad \quad k \in \mathcal{K}_{-1} \\
        \dfrac{dm_{{K}}^{(n+1)}}{dt} & = -\lambda_{{K}}^{*(n)}(t) m_{{K}}^{(n+1)}
    \end{split}
\end{equation}
similarly. Using the solution, find $\bar{p}^{(n+1)}(t)$:
\begin{equation}
    \bar{p}^{(n+1)}(t) = \dfrac{\sum_{k=1}^{{K}}m_{k}^{(n+1)}(t) p_{k}^{*(n)}(t)}{\sum_{k=1}^{{K}}m_k^{(n+1)}}.
\end{equation}
If $(\bar{p}^{(n+1)}(t),\ m^{(n+1)}(t))$ and $(\bar{p}^{(n)}(t),  \ m^{(n)}(t))$ are sufficiently close to each other, we terminate. Otherwise, we go back to Step 1 with the newly computed $(\bar{p}^{(n+1)}(t),\ m^{(n+1)}(t))$. \\
\subsubsection*{Remarks on the algorithm}
\begin{itemize}
    \item As a proxy for distance between $(\bar{p}^{(1)}(t), m^{(1)}(t))$ and $(\bar{p}^{(2)}(t), m^{(2)}(t))$, we take
    \begin{equation}
    \begin{split}
        L_{\bar{p}}(\bar{p}^{(1)}, \bar{p}^{(2)}) &:= \int_{0}^{T} \big|\bar{p}^{(1)}(t) - \bar{p}^{(2)}(t)\big|^2 dt, \\
        L_{m}(m^{(1)}, m^{(2)}) &:= \sum_{k=0}^{{K}} \int_{0}^{T} \big|m^{(1)}_{k}(t) - m^{(2)}_{k}(t)\big|^2  dt, \\
    \end{split}
    \end{equation}
\begin{equation}
    L((\bar{p}^{(1)}, m^{(1)}), (\bar{p}^{(2)}, m^{(2)})) := L_{\bar{p}}(\bar{p}^{(1)}, \bar{p}^{(2)}) + L_{m}(m^{(1)}, m^{(2)}).
\end{equation}
\item We calculate the time integrals for the error term using the same partition as we had while solving HJB and Kolmogorov equations.
\item We terminate the algorithm when the error $L((\bar{p}^{(1)}, m^{(1)}), (\bar{p}^{(2)}, m^{(2)}))$ defined above becomes smaller than $tol = 10^{-6}$.
\end{itemize}

\subsection{Numerical Stability and Convergence of the Algorithm}
To test the numerical stability of the algorithm and the uniqueness of the equilibrium point, we simulate a sample system with $N$-many random initial guesses and check the corresponding value functions $\{ V_k^{(1)}(t), \dots, V^{(N)}_k(t)\}$, price profiles $\{ \bar{p}_k^{(1)}(t), \dots, \bar{p}^{(N)}_k(t)\}$   and inventory distributions $\{ m_k^{(1)}(t), \dots, m^{(N)}_k(t)\}$. \\
For each quantity, for each time point $t$, and for each level $k > 0$, we calculate the \emph{``coefficient of variation"} (standard deviation / mean) from $N$ samples and finally, take the supremum of this quantity over the entire period, among all the levels. \\
Mathematically, for a generic quantity $\{q_k(t)\}_{k = 1, t = 0}^{{K}, T}$, we use the following formula:
\begin{equation}
    \text{CV}_q = \sup_{k\in\mathcal{K}, t\in[0, T]} \frac{\sqrt{\frac{1}{N} \Big [\sum_{n=1}^{N} \Big (q^{(n)}_k(t) - \dfrac{1}{N} \sum_{m = 1}^{N} q^{(m)}_k(t)\Big )^2 \Big ]}}{\frac{1}{N} \sum_{n = 1}^{N} q^{(n)}_k(t)}.
\end{equation}
For the values $N = 10$, $T = 200$, $r = 0.04$, $\epsilon = 0.4$, we find 
\begin{equation}
         \text{CV}_{V} \approx 4 \times 10^{-11},\ \  \text{CV}_{\bar{p}} \approx 5 \times 10^{-12},\ \ \text{CV}_{m} \approx 2 \times 10^{-8}.
\end{equation}
Based on our observations, the convergence is typically achieved in 6 steps, and the proportional decrease in error remains consistent after the first step. Notably, the proportionality constant changes with $\epsilon$, with smaller values of $\epsilon$ leading to faster convergence. This aligns with expectations, as for $\epsilon = 0$, convergence occurs in a single iteration due to the independence of the HJB solution from $m(t), \ \bar{p}(t)$. Furthermore, our experiments indicate that the algorithm requires roughly the same number of iterations to converge as $K$ increases, while the overall computation time grows quadratically in $K$.

\subsection{Algorithm Results}
We set ${K} = 100,\  T = 200$ days, $r = 0.04$ per day, with time discretization $\Delta t=T/1000$. For the initial distribution, we give a ``bimodal-like" distribution, comprising two groups of firms uniformly distributed over $k\in\{20,\dots,24\}\cup\{50,\dots,54\}$, that is, a low-inventory and a high-inventory group.
For our initial condition on $m$, we get the following graphs which shows the concavity of $V$ in $k$ at the starting time $t=0$, by demonstrating diminishing marginal returns principle in action: We see that incremental gains from additional inventory reduce as inventory level increases, due to saturation and increased urgency to sell before expiration.

\begin{figure}[htbp]
    \centering
    \includegraphics[width=0.5\textwidth]{ 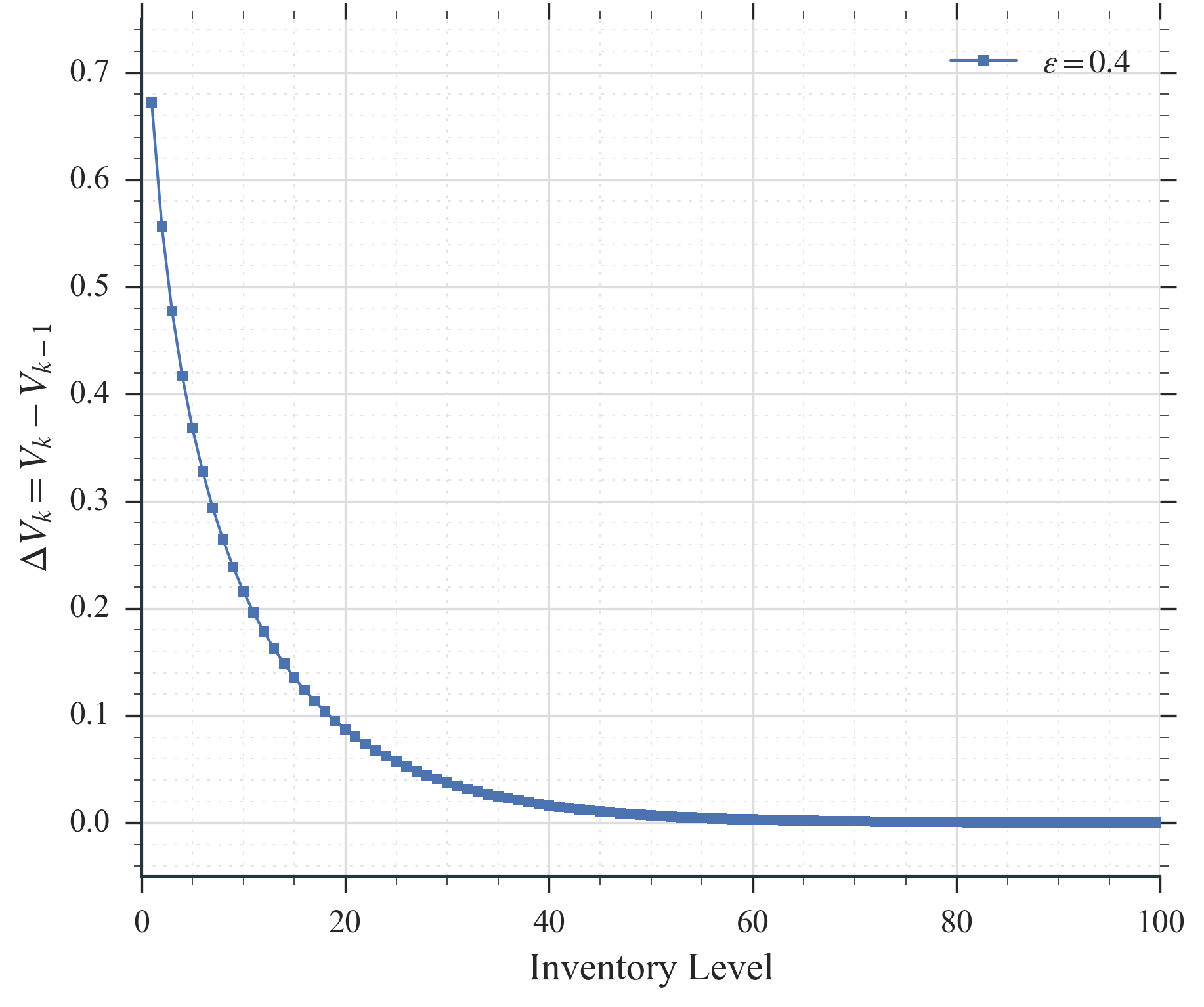}
    \caption{Concavity of $V$ in k}
    \label{fig:concavity}

\end{figure}

\begin{figure}[htbp]
    \centering
    \includegraphics[width=1\textwidth]{ 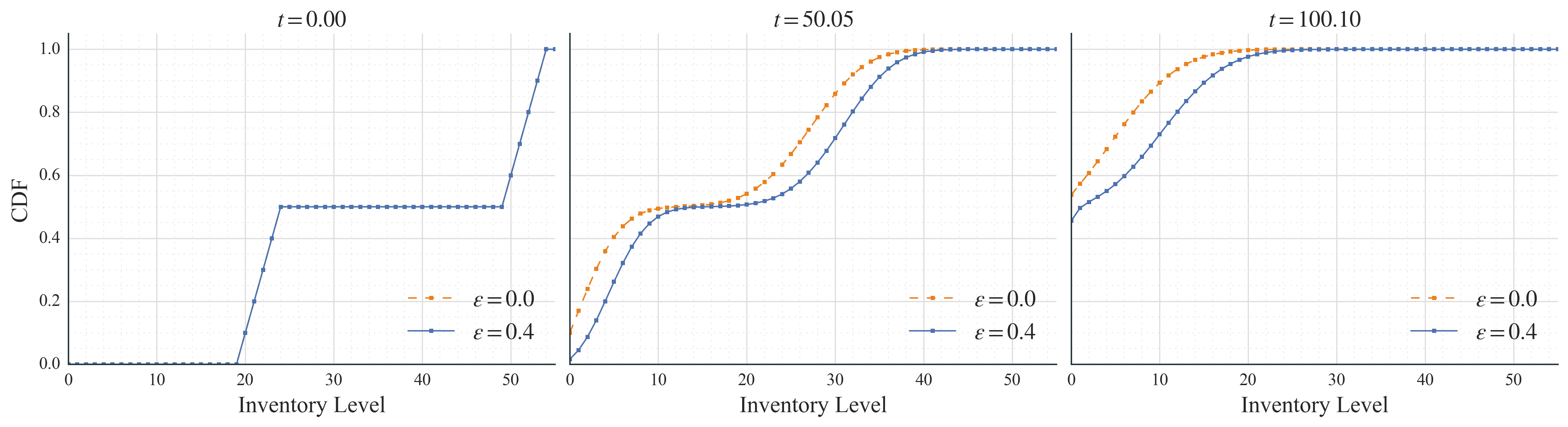}
    \caption{Evolution of CDF of $m(t)$ at given times for varying $\epsilon$. }
    \label{fig:m-evolution}
\end{figure}

In Figure \ref{fig:m-evolution}, we plot the cumulative distribution function (CDF) of $m$ at three time instances over inventory level $k$, for $\epsilon \in \{0,0.4\} $. The plots show that with increased competition, the market is slower, that is, the firms sell at a slower rate (evidenced by lagging densities when $\epsilon = 0.4$). As firms get more competitive, we get a more stable inventory flow. This suggests that intensified competition leads firms to behave cautiously, pricing tickets more conservatively to avoid rapid inventory depletion and subsequent competitive disadvantage. Hence, increased competition fosters more stable inventory dynamics, ensuring firms maintain availability throughout the selling period.

\begin{figure}[htbp]
    \centering
    \includegraphics[width=1\textwidth]{ 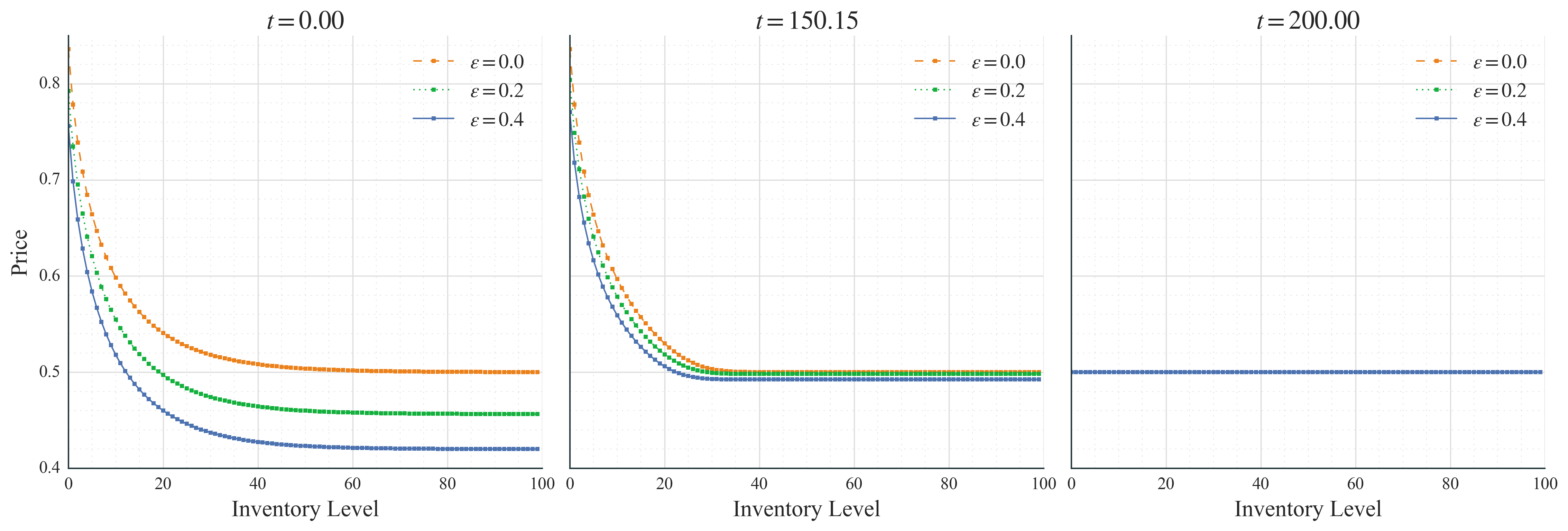}
    \caption{Evolution of optimal price $p^*(t)$ at given times for varying $\epsilon$.}
    \label{fig:pstar-over-time}
\end{figure}

Figure \ref{fig:pstar-over-time} reinforces this interpretation by displaying optimal price strategies under varying competition levels. Increased competition results in uniformly lower optimal prices across inventory levels, confirming competitive pressures incentivize price reductions to attract customers. Furthermore, higher inventory levels consistently translate to lower prices, indicating urgency in inventory liquidation as the time horizon shortens. This strategic pricing behavior highlights the necessity of balancing revenue generation with the risk of unsold inventory.

\begin{figure}[H]
    \centering
    %
    \begin{subfigure}[b]{0.47\textwidth} %
        \centering
        \includegraphics[width=\textwidth]{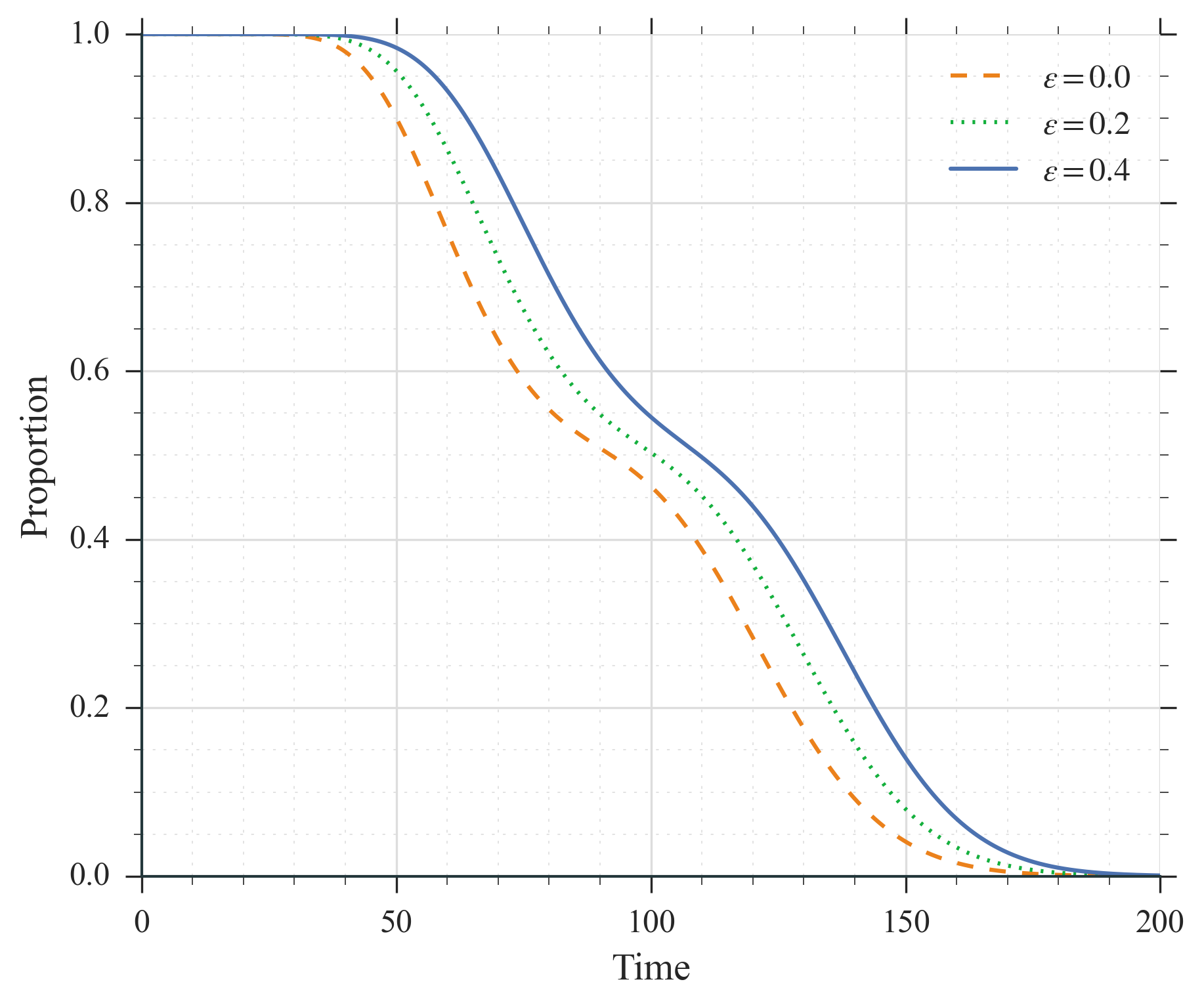}
        \caption{Ratio of active firms $\eta(t)$ for different $\epsilon$.}
        \label{fig:eta-vs-time}
    \end{subfigure}
    \hfill
    %
    \begin{subfigure}[b]{0.47\textwidth} %
        \centering
        \includegraphics[width=\textwidth]{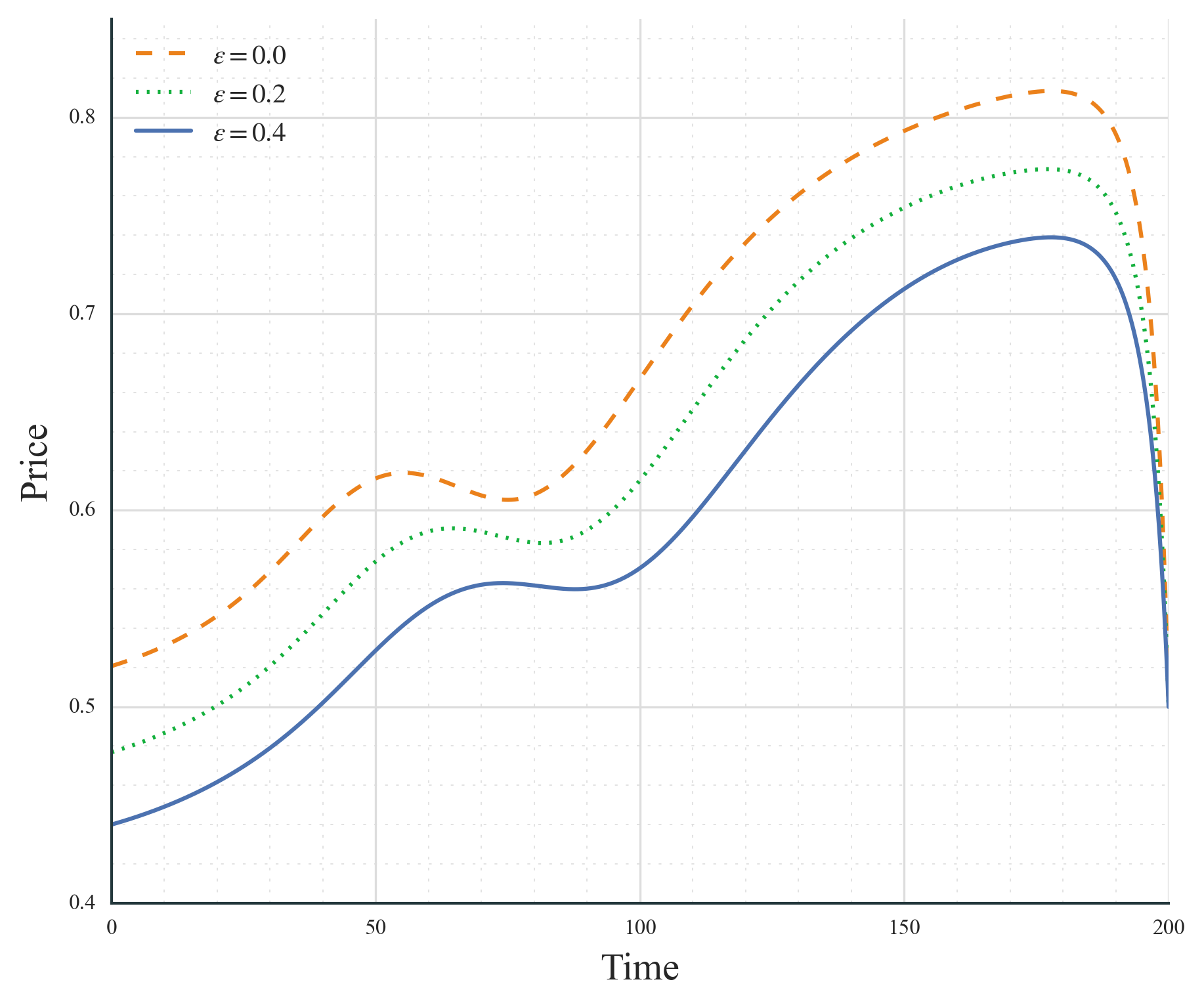}
        \caption{Average price $\bar{p}(t)$ for different $\epsilon$.}
        \label{fig:pbar-diff-eps}
    \end{subfigure}
    \caption{Comparisons for varying $\epsilon$ values.}
    \label{fig:eta-and-pbar}
\end{figure}
Figure \ref{fig:eta-and-pbar} displays the time-evolutions of two important market quantities, proportion of active firms $\eta(t)$ and average price in the market $\bar{p}(t)$ for varying $\epsilon$. The slower rate of selling with higher competition levels is also inferred from the proportion of active firms Figure \ref{fig:eta-vs-time}, as an increased proportion of firms are actively selling inventory over extended periods. In Figure \ref{fig:pbar-diff-eps}, as expected with more competition, average market price decreases, benefiting consumers. These two effects together underscore competition's role in prolonging market participation, allowing consumers greater purchasing flexibility.

In particular, the inflection points of the $\bar{p}(t)$ curve in Figure \ref{fig:pbar-diff-eps} correspond to key moments in the evolution of firm density, enhancing the interpretability of the results. The first inflection occurs when the low-inventory group begins to sell out; the second appears as this group approaches complete exit from the market; and the third inflection marks the start of market exit by the high-inventory group.

\begin{figure}[htbp]
    \centering
    \includegraphics[width=1\textwidth]{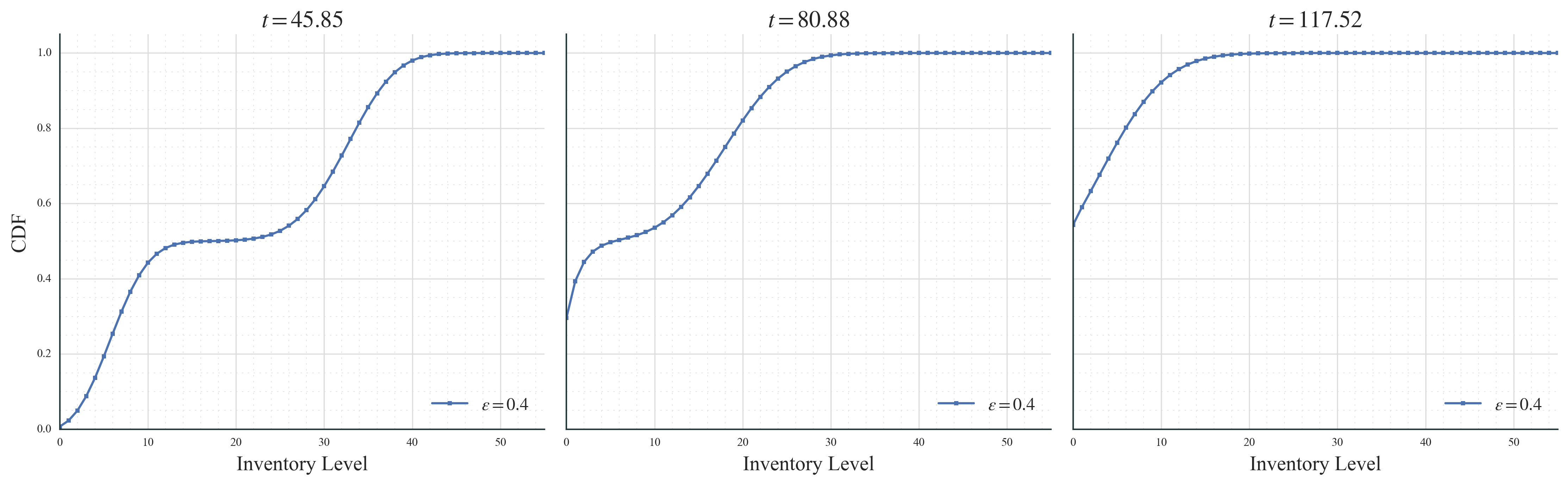}
    \caption{Density $m(t)$ at inflection points}
    \label{fig:inflection}
\end{figure}

\subsection{Qualitative Comparison with Data}

We intend our model to be viewed as an initial framework that demonstrates how dynamic competition and stochastic intensity control can be captured analytically in a mean-field setting. Much like \citep{Gallego1994}, we do not attempt calibration to data in this first work, but rather we establish the structural and theoretical foundations of this new class of models from which empirically relevant extensions can evolve. Detailed calibration to real-world airline data would require extensive and granular datasets that are not readily available (as effective demand intensities ($\lambda$ and related parameters) differ across routes), as well as the additional market functionalities (e.g., cancellations, refunds, customer learning) that go beyond the scope of the current paper.

\begin{figure}[htbp]
    \centering
    \includegraphics[width=0.8\textwidth]{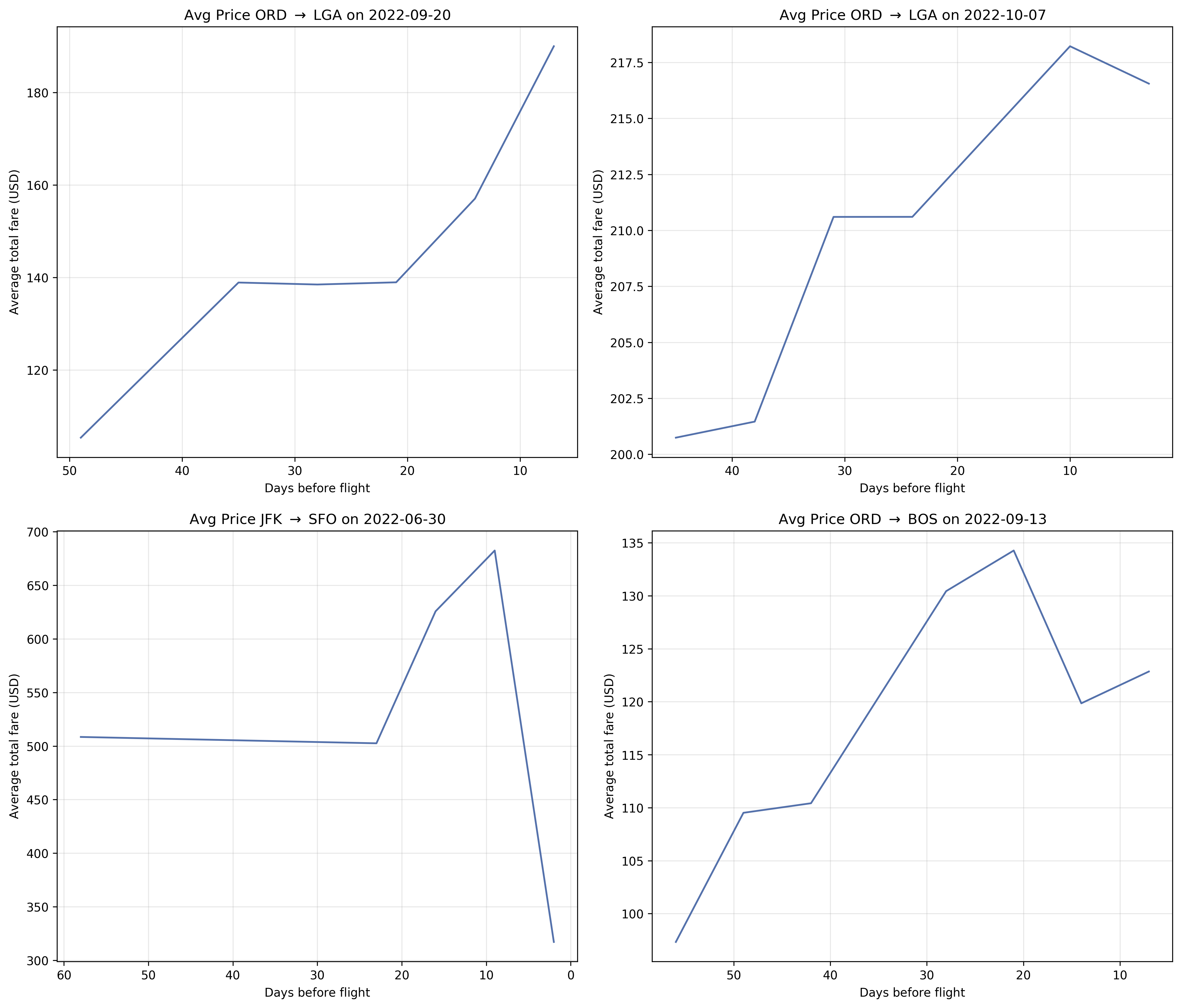}
    \caption{A survey of average flight prices across various dates and routes}
    \label{fig:survey-quad}
\end{figure}

Figure \ref{fig:survey-quad}, using 2022 airfare data \citep{kaggle}, illustrates broad qualitative consistency of our model with the general temporal variation of prices, leading up to the departure date. Across four sample routes shown, we see a pattern of increasing fares as the flight gets closer, with stagnant periods at intermediate times. However, as observed in the latter three plots with varying degrees of prominence, the prices may experience a drop shortly before the departure, consistent with our basic model (see Figure \ref{fig:pbar-diff-eps}). This phenomenon (and other similar features) are clearly quite dependent on the routes, and therefore our model requires careful adaptation and calibration to individual routes.

\section{Conclusion}\label{sec:concs}
In this paper, we introduced a novel existence and uniqueness result for a finite-state mean field game of stochastic intensity control, using differential equations methods. Our approach aims to provide a new perspective on analyzing such games in revenue management, and also to contribute additional tools for studying state mean field interactions in discrete-state settings.

We note that our ticket-selling model can be generalized to encompass a variety of situations seen in the real markets. One of these situations is the option of reselling by the customers. In this case, the price for reselling should also be another control parameter. This can be readily applied in our framework. Another example is the firm overselling. In this instance, firms have the option to oversell their tickets at the cost of paying high penalties at the end of the trading period if they oversell by that time. This example is seen in the airline companies for instance. Our model can be used to understand whether this scheme leads to economically good results for the customers similar to the firm competition.  In the computations for the overselling setup, the value functions of the following form can be considered: 
\begin{equation}\label{eq:discussion}
V_{k}(t) = \sup_{{p_s \geq 0}} \mathbbm{E} \Big \{  \int_{t}^{T} e^{-r(s-t)} p_s \mathbbm{1}_{\{X_s > 0\} } dN_s^\lambda + e^{-r(T-t)}\Psi(X_T)\Big | X_t = k\Big \},
\end{equation}
where $\Psi(X_T)$ is a penalty function for overselling, where $\Psi(x) = 0$ for $x \ge 0$, further tailored to suit the desired market model.

We can extend number of tickets being sold at every arrival by considering $N$ as a marked point process. In realistic situations, groups of people traveling together, such as families, would buy their tickets together at the same time, if there is sufficient inventory. Therefore, in our model, each arrival can be associated with a random number of sales, such that the HJB equation becomes
\begin{equation*}
    \dfrac{dV_k}{dt} - rV_k + \sup_{p\geq 0} \Big\{ \sum_{j=1}^{J \land k}\lambda_k\pi_j \big[jp-(V_k-V_{k-j})\big]\Big\} = 0,
\end{equation*}
and Kolmogorov equations can be written as
\begin{align*}
    \dfrac{dm_k}{dt} -\Big(\sum_{j=1}^{J \land K-k}\pi_j\lambda_{k+j}m_{k+j} - \lambda_k m_k\Big)= 0,
\end{align*}
subject to same initial and terminal conditions as in \eqref{eq:3.6}.

As we showed in Appendix \ref{app:4}, the existence and uniqueness of a fixed point can be proven for a sufficiently small interaction parameter $\epsilon$. Our proof relied upon topological fixed point theorems and an energy identity. Due to the complex dependence of the optimal controls on the input measure, it is not trivial to generalize the intensity function, as noted in the proof step where we rely on crucial properties of the base-case ($\epsilon = 0$) solution. Therefore, we restrict attention to the  piecewise linear intensity function in Section~\ref{sec:3.0}. More general intensity models, such as logit, probit, or isoelastic \citep[Chapter 6]{vives} can also be considered on a case-by-case basis.

We aim to extend our model to incorporate cancellations and refund mechanisms, as experienced in real world due to catastrophes or disruptions such as pandemics, weather anomalies, natural disasters, especially relevant for airline or hospitality industries. Further analysis can also be conducted in the presence of abundant data, both on its own using techniques such as reinforcement learning to approach this problem, and in addition to our method by calibrating the model parameters. 

Comparing the discrete-space proofs to their continuous-space analogs in the literature, the uniqueness proof would follow the same structure in a continuous-space mean field game, relying similarly on an energy identity and monotonicity arguments. The main difference arises in the existence proof: while our discrete-state setting allows us to use compactness arguments in finite dimensions, in the continuous-space case one would need to work with the Wasserstein distance on probability measures and establish additional regularity results.

Finally, for our future work, we plan to work on a proof that as the number of agents $n$ increases, we achieve an $\varepsilon_n$-Nash equilibrium that converges to the mean field equilibrium with $\varepsilon_n \to 0$. There is a rich literature on this subject for mean field games with Brownian motion, and in the most general circumstances, one achieves $\varepsilon_n = O(1/n)$ ( \citep{delarueramanan}, \citep{laurieretangpi}, \citep{nourian}). However, these theorems are not readily applicable to our system, and further studies need to be done on the subject. 

\appendix

\section{Brief Review and Limit of the Finite-Player Case} \label{app:3}

We show that the consumer preferences in our model can be derived from Quasilinear Quadratic Utility Model (QQUM). Letting $q = (q_i)_{i=1}^N$, a suitable utility function is
\begin{equation}
U(q) = \sum_{i=1}^{N} q_i - \frac{1}{2} \left( \sum_{i=1}^{N} q_i^2 + \frac{\epsilon}{N - 1} \sum_{j \ne i} q_i q_j \right).
\end{equation}
The first-order condition of the utility maximization problem $\displaystyle\max_{q \in \mathbbm{R}_+^N} U(q) - pq$ for a representative consumer yields the inverse demand functions: 
\begin{equation}
    q_i = D_i(p) = 
\begin{cases}
a_n - b_n p_i + c_n \bar{p}_i^n, & i = 1, 2, \cdots, n \\
0, & i = n+1, \cdots, N
\end{cases}
\end{equation}
where $n$-many firms receive positive demands, and thus, are active and produce positive quantities, and the rest of the firms are inactive \citep[Chapter 6 for more details]{vives}.
These constants can be calculated as 
\begin{equation}\label{eq:a3}
    \begin{split}
        &a_n = \frac{1}{1+\epsilon\frac{n-1}{N-1}}, \quad b_n = \frac{1+\epsilon\frac{n-2}{N-1}}{(1+\epsilon\frac{n-1}{N-1})(1-\frac{\epsilon}{N-1})}, \quad c_n=\frac{\epsilon\frac{n-1}{N-1}}{(1+\epsilon\frac{n-1}{N-1})(1-\frac{\epsilon}{N-1})}.
    \end{split}
\end{equation}
In the context of our paper, for very large ($n, \ N$), $\frac{n-1}{N-1} \to \eta(m)$ (the proportion of active firms), and $\frac{1}{N-1}\to0$. Consequently, the equations in \eqref{eq:a3} become
\begin{equation}
    \begin{split}
        a(m)=\frac{1}{1+\epsilon\eta(m)},\qquad b(m) = 1, \qquad c(m) = \frac{\epsilon\eta(m)}{1+\epsilon\eta(m)},
    \end{split}
\end{equation}
concluding the justification about our demand arrival process rate of $$\lambda(p, \overline{p},m) = (a(m) + c(m)\overline{p} - p)^{+}.$$

\section{Proof of Proposition \ref{prop:2}} \label{app:2}

\subsection{The inductive setup for Picard-Lindelöf theorem}
As the HJB equations become independent from the Kolmogorov equations when $\epsilon =0$, we shall start by first solving for $V_k^{(0)}$'s, proceeding by induction on $k$. Remember that the HJB equations are given by
\begin{equation} \label{eq:43}
    \dfrac{dV_k^{(0)}}{dt} - rV_k^{(0)} + \dfrac{1}{4} ( 1 - \Delta V_k^{(0)} )^2 = 0.
\end{equation}
Suppose, for an induction argument, that the bounds in Proposition \ref{prop:2} hold for some $k$. From \eqref{eq:43}, we get
\begin{equation} \label{eq:46}
    \dfrac{d\Delta V_{k+1}^{(0)} }{dt} + \dfrac{1}{4} (A_k - \Delta V_{k+1}^{(0)} ) (B_k - \Delta V_{k+1}^{(0)} ) = 0,
\end{equation}
where we defined 
\begin{equation*} \label{eq:45}
    \begin{split}
        & A_k(t) = 1 + 2r - \sqrt{(1 + 2r)^2 - 1 + (1 - \Delta V_{k}^{(0)}(t))^2} \\
        & B_k(t) = 1 + 2r + \sqrt{(1 + 2r)^2 - 1 + (1 - \Delta V_{k}^{(0)}(t))^2}.
    \end{split}
\end{equation*}

Now we argue that a unique solution to \eqref{eq:46} exists in $[0, T]$ and further it satisfies $0 \leq \Delta V_{k+1}^{(0)}(t) < A_k(t)$ for all $t\in [0, T]$. Introducing a new variable
\begin{equation*} \label{eq:47}
    z(t) = \dfrac{1}{A_k(t) - \Delta V_{k+1}^{(0)}(t)},
\end{equation*}
and $0 < z(T) \leq z(t) < \infty$ for all $t\in [0, T]$. Notice that for $z$, \eqref{eq:46} can be written as 
\begin{equation} \label{eq:48}
    \dfrac{dz}{dt} + \dfrac{1}{4} (B_k - A_k) z + A_k' z^2 + \dfrac{1}{4} = 0. 
\end{equation}

\subsection[Inductive step: existence and uniqueness of ΔVk(0)]{Inductive step, proving existence and uniqueness of $\Delta V_k^{(0)}$}
By Picard-Lindelöf theorem, it suffices to show that the solution does not blow up for any $t_0 \in [0, T]$. Therefore, we will now bound the solution from above and below.  
We can rewrite \eqref{eq:48} as
\begin{equation} \label{eq:411}
    z' + 2\alpha z - \beta z^2 + \gamma = 0, \nonumber
\end{equation}
where $\alpha(t), \beta(t), \gamma(t)$ are non-negative functions following from $B_k(t)-A_k(t)\geq 0$ and $A_k'(t)\leq 0$. Change of variables $\tilde{z}(t) = z(t)e^{\int_{0}^{t}\alpha(s) ds}$ yields 
\begin{equation} \label{eq:412}
    \tilde{z}' - \tilde{\beta} \tilde{z}^2 + \tilde{\gamma} = 0 \nonumber
\end{equation}
where $\tilde{\beta}(t) = \beta(t) e^{-\int_{0}^{t}\alpha(s) ds} \geq 0$ and $\tilde{\gamma}(t) = \gamma(t)e^{\int_{0}^{t}\alpha(s) ds} \geq 0$. Then the upper bound for $ \tilde{z}(t)$ can be shown by
\begin{equation} \label{eq:412a}
    \tilde{z}(t) \leq \tilde{z}(T) + \int_{0}^{T} \tilde{\gamma}(s) ds \nonumber
\end{equation}

Conversely for the lower bound, inspecting $({1/\tilde{z}})^{'} + \tilde{\beta} - \tilde{\gamma}/{\tilde{z}^2} = 0$ allows us to show that 
\begin{equation} \label{eq:417}
    \tilde{z}(t) \geq \Big ( \dfrac{1}{\tilde{z}(T)} + \int_{0}^{T} \tilde{\beta}(s)ds \Big)^{-1} \nonumber
\end{equation}

This concludes the proof that a solution $\Delta V_{k+1}^{(0)}$ in $[0, T]$ exists and is unique. Further for all $t\in[0, T]$
\begin{equation} \label{eq:418}
    \begin{split}
        0 \leq \Delta V_{k+1}^{(0)}(t) &\leq A_k(t)\leq 1 + 2r - 2 \sqrt{r^2 + r} = 1-2\delta\nonumber
    \end{split}
\end{equation}
and $\dfrac{d \Delta V_{k+1}^{(0)}}{dt} (t) \leq 0$, which is observed readily from \eqref{eq:48} and the fact that $B_k \geq A_k$. This concludes the inductive step. Now let us prove it for the base case. 

\subsection{Base case, proving existence and uniqueness of \texorpdfstring{$V_k^{(0)}$}{Vk(0)}}
In this case, we have 
\begin{equation} \label{eq:419}
    \dfrac{d \Delta V_1^{(0)}}{dt} + \dfrac{1}{4} (A_1 - \Delta V_1^{(0)})(B_1 - \Delta V_1^{(0)}) = 0, \nonumber
\end{equation}
where  
\begin{equation} \label{eq:420}
\begin{split}
    & A_1 = 1 + 2r - 2\sqrt{r^2 + r}\\
    & B_1 = 1 + 2r + 2\sqrt{r^2 + r}. \nonumber
\end{split}
\end{equation}
Since $A_1$ and $B_1$ have no dependence on $t$, the solution has an explicit form. 
\begin{equation} \label{eq:421}
    \Delta V_1^{(0)}(t) = A_1 \dfrac{1}{1 + \dfrac{B_1 - A_1}{A_1} (e^{\sqrt{r^2 + r} (T - t)}-1)}.  \nonumber
\end{equation}
It is seen that $\Delta V_1^{(0)}$ satisfies the appropriate conditions. This concludes the proof that there exists a unique solution to $V_k^{(0)}$ in $[0, T]$. 

\subsection{Proving existence and uniqueness of \texorpdfstring{$m_k^{(0)}$}{mk(0)}}
With the results for $V_k^{(0)}$ in mind, recall that the Kolmogorov equations are
\begin{equation} \label{eq:422}
\begin{split}
    &\dfrac{dm_k^{(0)}}{dt} - (\lambda_{k+1}^{(0)} m_{k+1}^{(0)} - \lambda_{k}^{(0)}m_k^{(0)})  = 0, \quad k \in \mathcal{K}_{-1} \\
    & \dfrac{dm_K^{(0)}}{dt} + \lambda_{K}^{(0)}m_K^{(0)}  = 0, \nonumber
\end{split}
\end{equation} 
with the initial condition $m_k^{(0)}(0) = M_k$, satisfying $0 \leq M_k \leq 1, \forall k \in \mathcal{K}$ and $\sum_{k=0}^K M_k = 1$.

Consequently, we have the following solutions
\begin{equation} \label{eq:423}
    \begin{split}
        & m_K^{(0)}(t) = M_K\exp {\Big ( - \int_{0}^{t} \lambda_K^{(0)} (s) ds \Big )} \\
        & m_k^{(0)} (t) = \exp {\Big ( - \int_{0}^{t} \lambda_k^{(0)} (s) ds \Big )} \Bigg[ M_k +\int_{0}^{t} \lambda_{k+1}^{(0)}(s)m_{k+1}^{(0)}(s) \exp{\Big ( \int_{0}^{s} \lambda_{k}^{(0)}(u)du \Big ) ds}\Bigg], \quad k\in\mathcal{K}_{-1}.
    \end{split}
\end{equation}

Since $M_K\geq0$ and all the other terms are nonnegative in \eqref{eq:423}, by induction, we can see that $\forall k \in \mathcal{K},\ \forall t \in [0,T], \ m_k^{(0)}(t) \geq 0$. 

Further, for $\eta^{(0)}(t;m^{(0)})= \sum_{k=1}^{K}m_k^{(0)}(t)$, we have 
\begin{equation} \label{eq:425}
    \dfrac{d\eta^{(0)}}{dt} = - \lambda_1^{(0)} m_1^{(0)} \nonumber,
\end{equation}
with $\eta^{(0)}(0;m^{(0)}) = 1$ (assuming without loss of generality that no firms start with zero inventory).  Therefore $\eta^{(0)}(t;m^{(0)}) \leq 1$ for $t\in [0, T]$, and this implies that for all $t\in [0, T]$
\begin{equation} \label{eq:426}
    0 \leq m_k^{(0)} (t) \leq 1, \quad k \in\mathcal{K} \nonumber.
\end{equation}
This concludes the proof.\qed


\section{Proofs of Existence and Uniqueness} \label{app:4}
\subsection{Proof of Theorem \ref{thm:2.2}}
Let us denote the monopolist solution as $(V^{(0)}, m^{(0)})$, and define the first-order deviations $V^{(1)} = \{ V_k^{(1)} \}_{k=1}^{K}$ and $m^{(1)} = \{ m_k^{(1)} \}_{k=1}^{K}$ by 
\begin{equation} \label{eq:3.7}
\begin{split}
    &V_k = V_k^{(0)} + \epsilon V_k^{(1)} \\
    &m_k = m_k^{(0)} + \epsilon m_k^{(1)},
\end{split}   
\end{equation}
where $V_k^{(1)}$ and $m_k^{(1)}$ will depend on $\epsilon$.
We will consider the vector space of continuous functions $$\mathcal{B} := \{(V_1^{(1)},\dots,V_K^{(1)}, m_1^{(1)},\dots,m_K^{(1)}),\ \ V_k^{(1)}, m_k^{(1)} \in C([0,T]), \ \ k\in\mathcal{K}\},$$ and the uniform norm
\begin{equation} \label{eq:218}
    \| (V^{(1)}, m^{(1)}) \|_{\mathcal{B}} = \sum_{k=1}^{K} \|V_k^{(1)}\|_{\infty} + \sum_{k=1}^{K} \|m_k^{(1)}\|_{\infty}.
\end{equation}
With this norm, $\mathcal{B}$ is a Banach space. Vector addition and scalar multiplication are defined elementwise and pointwise as usual. 

\subsubsection{Outline of the proof process}

By formulating two operators $R$ and $S$, we will transform our problem into finding a fixed point $(V^{(1)}$, $m^{(1)})$, such that
\begin{equation} \label{eq:217}
\begin{split}
    &V_k^{(1)} = R_k (V^{(1)}, m^{(1)}), \quad k \in \mathcal{K} \\
    &m_k^{(1)} = S_k(V^{(1)}, m^{(1)}), \quad k \in \mathcal{K}.
\end{split}
\end{equation}
To use the Schauder Fixed Point Theorem, we will first define the mapping $\Gamma$:
\begin{equation} \label{eq:219}
    \Gamma((V^{(1)}, m^{(1)})) = (R_1(V^{(1)}, m^{(1)}), ..., R_K(V^{(1)}, m^{(1)}), S_1(V^{(1)}, m^{(1)}), ... , S_K(V^{(1)}, m^{(1)})).
\end{equation}
Then, we will find the bounds $\{C_{k}^{V}, C_{k}^{m}\}_{k=1}^{K}$ such that the following closed and convex ellipsoid
\begin{equation} \label{eq:220}
    M := \{(V^{(1)}, m^{(1)}) \in \mathcal{B} : \|V_k^{(1)}\|_{\infty} \leq C_k^{V}, \ \|m_k^{(1)}\|_{\infty} \leq C_k^{m}, \ k \in \mathcal{K} \}
\end{equation}
is mapped into itself under $\Gamma$. Thereafter by proving that $\Gamma$ is continuous and $\Gamma(M)$ is precompact (i.e. its closure is compact), we will prove the existence of a fixed point.

In the proof, we will derive the operators and define the bounds $\{C_k^{V} \}_{k=1}^{K}$ and $\{C_k^{m} \}_{k=1}^{K}$ as functions solely dependent on the external parameters: $T, K, r$. Subsequently, we will prove that for a sufficiently small $\epsilon > 0$, the mapping $\Gamma$ maps $M$ into itself. Next, we will establish the continuity of $\Gamma$ with respect to the norm $\|\cdot\|_{\mathcal{B}}$. And finally, we will show that the family $\Gamma(M)$ of functions is equicontinuous and pointwise-bounded, which will allow us to conclude that $\Gamma(M)$ is precompact by the Arzela-Ascoli Theorem.

\subsubsection{Deriving the integral operators}

Our strategy is to linearize \eqref{eq:3.6} and solve it using Schauder's Fixed Point Theorem. To begin, let us establish the following notations:
\begin{equation}\label{eq:3.8a}
    \begin{split}
        &\eta^{(0)} =\sum_{k=1}^{K} m_k^{(0)}, \qquad\eta^{(1)} =\sum_{k=1}^{K} m_k^{(1)}, \\
        &\phi^{(0)} = \phi_0(V^{(0)}, m^{(0)}) = \frac{1}{2}\Big(\sum_{k=1}^{K} m_k^{(0)} \Delta V_k^{(0)} - \eta^{(0)} \Big ).
    \end{split}
\end{equation}
Aiming to satisfy $\phi_\epsilon(V^{(1)}, m^{(1)}) = \phi^{(0)} + \epsilon \phi^{(1)}$, we define the first-order perturbation term $\phi^{(1)}$ as 
\begin{equation} \label{eq:3.9}
    \phi^{(1)} = \dfrac{1}{2 + \epsilon (\eta^{(0)} + \epsilon \eta^{(1)})} \Big (- (\eta^{(0)} + \epsilon \eta^{(1)}) \phi^{(0)} +  \sum_{k=1}^{K} m_k^{(1)} \Delta V_k^{(0)}  + \sum_{k=1}^{K} m_k^{(0)} \Delta V_k^{(1)} - \eta^{(1)} \Big ) 
\end{equation}
Then the first equation in \eqref{eq:3.6} becomes
\begin{equation} \label{eq:3.10}
\begin{split}
    \dfrac{dV_k^{(1)}}{dt} - (\lambda_k^{(0)} + r) V_k^{(1)} + \lambda_k^{(0)} (\phi^{(0)} + V_{k-1}^{(1)}) + \epsilon \Lambda_{k}^{V} = 0,
\end{split}
\end{equation} 
where $\lambda_k^{(0)}$ and $\Lambda_k^{V}$ are defined in the following way:
\begin{equation} \label{eq:3.11}
\begin{split}
    & \lambda_k^{(0)} = \dfrac{1}{2} (1 - \Delta V_k^{(0)}) \\
    & \Lambda_k^{V} = \lambda_k^{(0)} \phi^{(1)} + \dfrac{1}{4} \Big (\phi^{(0)} - \Delta V^{(1)}_k + \epsilon \phi^{(1)} \Big )^2.
\end{split}
\end{equation}
Let us introduce the integral operator 
\begin{equation} \label{eq:3.12}
    R_k (V^{(1)}, m^{(1)}) (t) = \int_{t}^{T} e^{-\int_{t}^{s} (\lambda_k^{(0)}(s') + r) ds'} \Big ( \lambda_k^{(0)} (\phi^{(0)} + V_{k-1}^{(1)}) + \epsilon \Lambda_k^{V} \Big ) ds.
\end{equation}
With this definition, our problem is transformed into finding a fixed point $V^{(1)}$ under $R$, that is, for $ k \in \mathcal{K}$,  $V_k^{(1)} = R_k (V^{(1)}, m^{(1)})$.

Similarly, let us derive an integral operator for $m^{(1)}$. We start by defining
\begin{equation} \label{eq:3.14}
    \lambda_k^{(1)} = \dfrac{1}{2} \Big ( \phi^{(0)} - \Delta V^{(1)}_k + \epsilon \phi^{(1)}\Big ),
\end{equation}
the first-order deviation from $\lambda_k$ such that $\lambda_k \equiv \dfrac{1}{2} (1 - \Delta V_k + \epsilon \phi_\epsilon) = \lambda_k^{(0)} + \epsilon \lambda_k^{(1)} $. We observe that positivity of $\lambda_k$ is preserved for small enough $\epsilon > 0$.

We further define
\begin{equation} \label{eq:3:13}
     \Lambda_k^{m} = m_{k+1}^{(1)} \lambda^{(1)}_{k+1} - m_{k}^{(1)} \lambda^{(1)}_{k},\text{ for } k \in \mathcal{K}_{-1} \text{ and }\Lambda_K^{m}  = - m_{K}^{(1)} \lambda^{(1)}_{K}
\end{equation}
to simplify the Kolmogorov equations in \eqref{eq:3.6} and arrive at the integral operator:
\begin{equation} \label{eq:216}
    \begin{split}
        & S_k(V^{(1)}, m^{(1)}) (t) = \int_{0}^{t} e^{- \int_{s}^{t} \lambda_k^{(0)}(s')ds'} \Big ( m_{k+1}^{(1)}\lambda_{k+1}^{(0)} +  m_{k+1}^{(0)} \lambda^{(1)}_{k+1} - m_k^{(0)} \lambda^{(1)}_k + \epsilon \Lambda_{k}^m\Big ) ds\ \quad k\in \mathcal{K}_{-1}\\
        & S_K(V^{(1)}, m^{(1)})(t) = \int_{0}^{t} e^{- \int_{s}^{t} \lambda_k^{(0)}(s')ds'} \Big ( - m_K^{(0)} \lambda^{(1)}_K + \epsilon \Lambda_{K}^m\Big ) ds.
    \end{split}
\end{equation}

\subsubsection{Bounding the operators from above}

In order to show the results of existence, we determine that we can set $\{C_k^V \}_{k=1}^K$ and $\{C_k^m \}_{k=1}^K$ recursively as follows:
\begin{equation} \label{eq:222}
    \begin{split}
        C_1^V = 2K\quad\ \ \  \ \ \ \ \ \  &C_k^V = 2 C_{k-1}^V + 2K \text{, } \qquad \quad \ k\in\{2,\dots,K\} \\
        C_K^m = TK2^{K+5}\quad\ \  &C_k^m = T C_{k+1}^m + TK2^{K+5}, \ \ \ k\in\mathcal{K}_{-1}.
    \end{split}
\end{equation}
We observe that for all $k \in \mathcal{K}$, we have $C_k^V \leq C_K^V$ and $C_k^m \leq C_0^m$. From the definition of $\phi^{(0)}$ \eqref{eq:3.8a}, we have
\begin{equation} \label{eq:229}
        \|\phi^{(0)}\|_{\infty} \leq \dfrac{1}{2} \Big ( \sum_{k=1}^{K} \|m_k^{(0)}\|_{\infty}\|\Delta V_k^{(0)}\|_{\infty} + \|\eta^{(0)}\|_{\infty}\Big ) \leq K.
\end{equation}
In order to bound $\phi^{(1)}$, we take $\epsilon > 0$ small enough such that
\begin{equation} \label{eq:228}
\begin{split}
     \|\phi^{(1)}\|_{\infty} & \leq \Bigg\|\dfrac{1}{2 + \epsilon((\eta^{(0)} + \epsilon \eta^{(1)}) )} \Bigg\| _{\infty} \Bigg (\|\eta^{(0)} + \epsilon \eta^{(1)}\|_{\infty} \|\phi^{(0)}\|_{\infty} +\\
     & \sum_{k=1}^{K} \|m_k^{(1)}\|_{\infty} \|\Delta V_k^{(0)}\|_{\infty} + \sum_{k=1}^{K} \|m_k^{(0)}\|_{\infty} \|\Delta V_k^{(1)}\|_{\infty} \\ &+ \epsilon \sum_{k=1}^{K} \|m_k^{(1)}\|_{\infty} \|\Delta V_k^{(1)}\|_{\infty} +\|\eta^{(1)}\|_{\infty} \Bigg) \\
     & \leq 2K C_0^m + 4K C_K^{V},
\end{split}
\end{equation}
where we used the facts that $0\leq \eta^{(0)}(t) \leq 1$ for all $t\in [0,T]$ and $ 0 \leq m_k^{(0)}(t) \leq 1$ for all $ k \in \mathcal{K}$, $t\in  [0,T]$ (proved in Appendix \ref{app:2}).
Similarly, for $\{ \lambda_k^{(1)} \}_{k=1}^{K}$ we get
\begin{equation} \label{eq:231}
\begin{split}
    \|\lambda_k^{(1)}\|_{\infty} \leq \frac{K}{2} + C_K^V + \dfrac{\epsilon}{2} \|\phi^{(1)}\|_{\infty}.
\end{split}
\end{equation}
Consequently, 
\begin{equation} \label{eq:232}
    \begin{split}
        \|S_k(V^{(1)}, m^{(1)})\|_{\infty} & \leq T \Big ( \|\lambda_{k+1}^{(0)}\|_{\infty} \|m_{k+1}^{(1)}\|_{\infty} + \|m_{k+1}^{(0)}\|_{\infty} \|\lambda^{(1)}_{k+1}\|_{\infty}+ \\
        &\qquad\qquad\qquad\|m_{k}^{(0)}\|_{\infty} \|\lambda^{(1)}_{k}\|_{\infty} + \epsilon \|\Lambda_k^m\|_{\infty} \Big ) \\
        & \leq \dfrac{T}{2} \Big ( C_{k+1}^{m} + 2K + 4 C_K^V \Big ) + \epsilon T \Big ( \|\phi^{(1)}\|_{\infty} + \|\Lambda_k^m\|_{\infty}\Big ),
    \end{split}
\end{equation}
where we used the bound $\|\lambda_k^{(0)}\|_{\infty} \leq 1/2$ (Appendix \ref{app:2}). As $\Lambda_k^m$ can be bounded by
$\|\Lambda_{k}^m\|_{\infty} \leq 2 C_0^m (C_K^V + \dfrac{\epsilon}{2} \|\phi^{(1)}\|_{\infty})$
, we can take $\epsilon  > 0$ small enough such that for all $k \in \mathcal{K}_{-1}$,
\begin{equation} \label{eq:234}
    \epsilon (\|\phi^{(1)}\|_{\infty} + \|\Lambda^{m}_k\|_{\infty} ) \leq \dfrac{1}{2} (C_{k+1}^m + 2K + 4C_K^V).
\end{equation}
Then we get, for $k \in \mathcal{K},$
\begin{equation} \label{eq:235}
\begin{split}
    \|S_k(V^{(1)}, m^{(1)})\|_{\infty} \leq T \Big ( C_{k+1}^{m} + 2K + 4 C_K^V \Big ) \leq C_k^{m},
\end{split}
\end{equation}
by observing that, for example, $\|\Lambda_{K}^m\|_{\infty}$ can be bounded above with the same constant we use for $\|\Lambda_{K-1}^m\|_{\infty}$ to conclude for the case when $k=K$. This is what we wanted. 

Similarly, for integral operators $\{R_k\}_{k=1}^K$ \eqref{eq:3.12} related to HJB equations, for sufficiently small $\epsilon > 0$ we have 
\begin{equation} \label{eq:238}
\begin{split}
    \|R_k(V^{(1)}, m^{(1)})\|_{\infty} & \leq \|\phi^{(0)} + V_{k-1}^{(1)}\|_{\infty} + \epsilon \Bigg\| \dfrac{1}{\lambda_{k}^{(0)}}\Bigg\|_{\infty} \|\Lambda_k^V\|_{\infty} \\
    & \leq \dfrac{C_{k}^V}{2} + \epsilon \Bigg\| \dfrac{1}{\lambda_{k}^{(0)}}\Bigg\|_{\infty}\|\Lambda_k^V\|_{\infty}\leq C_k^V,
\end{split}
\end{equation}
using the fact (proved in Appendix \ref{app:2}) that $\displaystyle \max_{k \in \mathcal{K}} \Bigg \| \dfrac{1}{\lambda_k^{(0)}} \Bigg \| _{\infty} = \frac{1}{\rho} < \infty$.

To find an explicit bound for $\epsilon$, note that for the case of $T<1$, by using explicit expressions for $\{C_k^V \}_{k=1}^K$ and $\{C_k^m \}_{k=1}^K$, we need $\epsilon < \dfrac{1 + KT}{K^2T2^{K+8}}$ to satisfy \eqref{eq:234}, and $\epsilon < \dfrac{\sqrt{r}}{{K^2T2^{K+8}} }$ to satisfy \eqref{eq:238}. Combining these, we get 
\begin{equation}
    \epsilon < C_1(K, T, r) = \dfrac{1}{K^2T2^{K+8}} \min \Big \{ 1 + KT, \sqrt{r}\Big \}.
\end{equation}
The case for $T>1$ can be handled similarly to calculate the expression for $C_1$.

\subsubsection{Applying Schauder Fixed Point Theorem}
Now, we have $\Gamma(M) \subset M$. The continuity of $\Gamma: M \to M$ is clear. For compactness, observe that the integrands in $R_k$ \eqref{eq:3.12} and $S_k$ \eqref{eq:216} can be uniformly bounded. Namely, we have a constant $C$ that depends only on the external parameters $T, K, r$ such that for all $(V^{(1)}, m^{(1)}) \in M$, for all $t,s\in [0,1]$ and for all $k \in \mathcal{K}$, we have 
\begin{equation}  \label{eq:247}
    \begin{split}
        &|R_k(V^{(1)}, m^{(1)})(t) - R_k(V^{(1)},m^{(1)})(s)| \leq C |t-s| \\
        & |S_k(V^{(1)}, m^{(1)})(t) - S_k(V^{(1)}, m^{(1)})(s)| \leq C |t-s|.
    \end{split}
\end{equation}
This shows that the family of functions $\Gamma(M)$ is uniformly equicontinuous. Since we are using a uniform norm on our Banach space, the family $\text{clo} (\Gamma(M))$ is also a uniformly equicontinuous family, and since $M$ is a closed set, we have $\text{clo} (\Gamma(M)) \subset M$. Hence $\text{clo} (\Gamma(M))$ is uniformly bounded. Consequently, by the $d$-dimensional Arzela-Ascoli Theorem, $\text{clo} (\Gamma(M))$ is a compact set. Then by the Schauder Fixed Point Theorem, we conclude that there exists $(V^{(1)}, m^{(1)}) \in M$ such that for all $t\in [0,T]$
\begin{equation} \label{eq:248}
    V_k^{(1)}(t) = \int_{t}^{T} e^{-\int_{t}^{s} (\lambda_k^{(0)} + r) ds'} \Big ( \lambda_k^{0} (\phi^{(0)} + V_{k-1}^{(1)}) + \epsilon \Lambda_k^{V} (V^{(1)}, m^{(1)})\Big ) ds, \quad k \in \mathcal{K}
\end{equation}
and 
\begin{equation} \label{eq:249}
    \begin{split}
        & m_k^{(1)} (t) = \int_{0}^{t} e^{- \int_{s}^{t} \lambda_k^{(0)}(s')ds'} \Big ( \lambda_{k+1}^{(0)}m_{k+1}^{(1)} - \Big ( m_{k+1}^{(0)} \lambda^{(1)}_{k+1} - m_k^{(0)} \lambda^{(1)}_k\Big ) + \epsilon \Lambda_{k}^m(V^{(1)}, m^{(1)})\Big ) ds, \ \  k \in \mathcal{K}_{-1}\\
        & m_K^{(1)}(t) = \int_{0}^{t} e^{- \int_{s}^{t} \lambda_k^{(0)}(s')ds'} \Big ( m_K^{(0)} \lambda^{(1)}_K + \epsilon \Lambda_{K}^m(V^{(1)}, m^{(1)})\Big ) ds.
    \end{split}
\end{equation} 

Differentiating both sides of\eqref{eq:248} and \eqref{eq:249} with respect to $t$, we conclude that 
\begin{equation} \label{eq:250}
\begin{split}
    V_k(t) := V_k^{(0)}(t) + \epsilon V_k^{(1)}(t), \quad &k \in \mathcal{K}\\
    m_k(t) := m_k^{(0)}(t) + \epsilon m_k^{(1)}(t), \quad &k \in \mathcal{K}
\end{split}
\end{equation}
is a solution to our set of coupled ordinary differential equations \eqref{eq:3.6}. \qed 

\subsection{Proof of Theorem \ref{thm:2.3}}
It is important to note that the Schauder Fixed Point Theorem does not guarantee uniqueness. In fact, establishing a priori uniqueness of a solution is not evident from the outset in the context of economic considerations. Nevertheless, when the interactions between the firms are sufficiently small, we can indeed prove uniqueness.
\subsubsection{Setup and result for when \texorpdfstring{$\epsilon = 0$}{epsilon=0}}
We shall use an energy identity to prove that the solution is unique. To begin, suppose that $(V^{[1]}, m^{[1]})$ and $(V^{[2]}, m^{[2]})$ are two solutions to the equations:
\begin{equation} \label{eq:31}
    \begin{split}
        & \dfrac{dV_k}{dt} - rV_k + H_k = 0, \quad k \in \mathcal{K} \\ 
        & \dfrac{dm_k}{dt} - (G_{k+1}m_{k+1} - G_k m_k) = 0, \quad k \in \mathcal{K}
    \end{split}
\end{equation}
where, for $k \in \mathcal{K}$, we defined
\begin{equation} \label{eq:32}
\begin{split}
    &G_k = \dfrac{1}{2} \Big ( 1 - \Delta V_k + \epsilon \phi_\epsilon \Big ) \quad H_k = G_k^2 = \dfrac{1}{4} \Big ( 1 - \Delta V_k + \epsilon \phi_\epsilon \Big )^2,
\end{split} 
\end{equation}
with the understanding that $G_{K+1} = 0$. For the sake of notational brevity, let us denote $\delta L := L^{[2]}- L^{[1]}$ for any generic quantity $L$. Applying $\delta$ to \eqref{eq:31}, multiplying the first equation by $\delta m_k$ and the second by $\delta V_k$, and adding them yields
\begin{equation} \label{eq:34}
\begin{split}
     2\dfrac{d \delta m_k \delta V_k}{dt} - r \delta m_k \delta V_k + \delta m_k \delta H_k - \delta V_k (\delta (G_{k+1}m_{k+1}) - \delta(G_k m_k)) = 0.
\end{split}
\end{equation}

Since the two solutions agree on the initial and terminal conditions, we have
\begin{equation} \label{eq:35}
    \delta V_k (T) = 0, \quad  \delta m_k (0) = 0,  \quad k \in \mathcal{K}.
\end{equation}
We integrate \eqref{eq:34} and sum over all $k \in \mathcal{K}$. Then, by discrete integration-by-parts, we get

\begin{equation} \label{eq:39}
    \int_{0}^{T} e^{-rt} \sum_{k=1}^{K} \Big \{ \delta m_k \delta H_k + \delta(\Delta V_k) \delta (G_k m_k) \Big \}dt = 0.
\end{equation}
We will now individually target the terms. To simplify this expression, let us introduce another notation: For any generic quantity $L$, denote $\overline{L} := (L^{[1]} + L^{[2]})/{2}$. We have the identity $\delta (L F) = \overline{F} \delta L + \overline{L} \delta F$ for another generic quantity $F$. Using the definitions for $H_k$ and $G_k$ yields $\forall k \in \mathcal{K},$
\begin{equation} \label{eq:312}
\begin{split}
    &\delta H_k = 2 \overline{G_k} \delta G_k  = -\dfrac{1}{2} \delta (\Delta V_k)\Big ( 1 - \overline{\Delta V_k} \Big ) + \epsilon Z_k \\
    &\delta (G_k m_k ) = \dfrac{1}{2} \Big ( \delta m_k (1 - \overline{\Delta V_k}) - \delta (\Delta V_k) \overline{m_k})\Big ) + \epsilon \zeta_k,
\end{split}
\end{equation}
where we defined
\begin{equation} \label{eq:313}
\begin{split}
    &Z_k = -\dfrac{1}{2} \Big ( \overline{\phi_\epsilon} \delta(\Delta V_k) - \delta \phi_\epsilon (1 - \overline{\Delta V_k})\Big ) + \dfrac{\epsilon\overline{\phi_\epsilon}}{2} \delta \phi_\epsilon, \quad \zeta_k = \dfrac{1}{2} \Big ( \overline{\phi_\epsilon} \delta m_k + \overline{m_k} \delta \phi_\epsilon \Big ).
\end{split}
\end{equation}
Then, \eqref{eq:39} becomes
\begin{equation} \label{eq:317}
    \dfrac{1}{2}\int_{0}^{T} e^{-rt} \sum_{k=1}^{K} \overline{m_k} (\delta(\Delta V_k))^2dt = \epsilon \int_{0}^{T}  e^{-rt}\sum_{k=1}^{K} \Big \{Z_k \delta m_k + \zeta_k \delta(\Delta V_k) \Big \} dt .
\end{equation}
Notice that for $\epsilon = 0$, we can conclude that $V^{[1]} = V^{[2]}$ by non-negativity of the left-hand-side integrand. Then, after applying $\delta$ to second line in \eqref{eq:31}:
\begin{equation}\label{eq:317a}    
\dfrac{d\delta m_k}{dt} - (\delta (G_{k+1}m_{k+1}) - \delta(G_k m_k)) = 0, \quad k \in \mathcal{K},
\end{equation}
we can observe that $m^{[1]} = m^{[2]}$ as well, when $\epsilon =0$. \\

\subsubsection{Uniqueness of \texorpdfstring{$V$}{V} in the general case}
For the general case when $\epsilon > 0$, we need to compute the expression on the right-hand side of \eqref{eq:317} explicitly. To begin, let us observe that 
\begin{equation} \label{eq:318}
    \begin{split}
        & Z_k \delta m_k + \zeta_k \delta(\Delta V_k) = \dfrac{\delta \phi_\epsilon}{2} \Big ( \overline{m_k} \delta(\Delta V_k) + \delta m_k (1 - \overline{\Delta V_k}) \Big ) + \dfrac{\epsilon \overline{\phi_\epsilon}}{2} \delta \phi_\epsilon \delta m_k.
    \end{split}
\end{equation}
Therefore, 
\begin{equation} \label{eq:319}
    \begin{split}
        \sum_{k=1}^{K} Z_k \delta m_k + \zeta_k \delta(\Delta V_k) = \dfrac{1}{2} \delta \phi_\epsilon \sum_{k=1}^{K}  \Big (\overline{m_k} \delta(\Delta V_k) + \delta m_k (1 - \overline{\Delta V_k}) \Big ) + \dfrac{\epsilon\overline{\phi_\epsilon}}{2}  \delta \phi_\epsilon \delta \eta.
    \end{split}
\end{equation}
Now let us find $\delta \phi_\epsilon$. Recall that we have 
\begin{equation} \label{eq:320}
    (2 + \epsilon \eta) \phi_\epsilon = \sum_{k=1}^{K} m_k (1 - \Delta V_k).
\end{equation}
Applying the difference operator $\delta$ on both sides, we have 
\begin{equation} \label{eq:321}
    (2 + \epsilon \overline{\eta}) \delta \phi_\epsilon + \epsilon \overline{\phi_\epsilon} \delta \eta = \sum_{k=1}^{K} \Big ( \overline{m_k} \delta (\Delta V_k) - \delta m_k (1 - \overline{\Delta V_k}) \Big )
\end{equation}
and thus,
\begin{equation} \label{eq:321a}
   \delta \phi_\epsilon = \dfrac{1}{ 2 + \epsilon \overline{\eta} } \Bigg[\sum_{k=1}^{K} \Big ( \overline{m_k} \delta (\Delta V_k) - \delta m_k (1 - \overline{\Delta V_k}) \Big ) -\epsilon \overline{\phi_\epsilon} \delta \eta\Bigg].
\end{equation}
Substituting $\delta\phi_\epsilon$ from \eqref{eq:321a} into the right-hand side of \eqref{eq:319} transforms the equation to
\begin{equation} \label{eq:323}
    \begin{split}
        & \sum_{k=1}^{K} Z_k \delta m_k + \zeta_k \delta(\Delta V_k) \\
        & = \dfrac{1}{2(2 + \epsilon \overline{\eta})} \Big ( \sum_{k=1}^{K} \overline{m_k} \delta(\Delta V_k) \Big )^2 - \dfrac{1}{2(2 + \epsilon \overline{\eta})} \Big [ \Big (\sum_{k=1}^{K} \delta m_k (1 - \overline{\Delta V_k})\Big) + \epsilon \overline{\phi_\epsilon} \delta \eta\Big]^2 \\
        & \leq \dfrac{K}{2(2 + \epsilon)} \sum_{k=1}^{K} \overline{m_k} (\delta(\Delta V_k))^2,
    \end{split}
\end{equation}
where in the last line we used that $\overline{m_k}(t) \in[0,1]$ for all $k\in \mathcal{K},\  t\in [0,T]$ and $\overline{\eta}(t) \in [0, 1]$ for all $t\in [0,T]$. Now combining (\ref{eq:317}) and (\ref{eq:323}), we obtain
\begin{equation} \label{eq:326}
    \dfrac{1}{2}\int_{0}^{T} e^{-rt} \sum_{k=1}^{K} \overline{m_k} (\delta(\Delta V_k))^2 dt \leq \dfrac{\epsilon K}{2(2 + \epsilon) } \int_{0}^{T} e^{-rt} \sum_{k=1}^{K} \overline{m_k} (\delta(\Delta V_k))^2dt.
\end{equation}
Consequently, for $K > 1$, and a sufficiently small positive $\epsilon < C_2 := \dfrac{2}{K-1}$, we get 
\begin{equation} \label{eq:327}
    \int_{0}^{T} e^{-rt} \sum_{k=1}^{K} \overline{m_k} (\delta(\Delta V_k))^2 dt = 0.
\end{equation}
Since $\overline{m_k}(t) > 0$ for $t\in (0, T]$, $ k \in \mathcal{K}$, we get $\delta(\Delta V_k(t))^2 = 0 \text{ for all } t\in (0, T]$. From continuity of $V_k(t)$ in $t$, we have $\delta(\Delta V_k)(0) = 0$ as well. Therefore 
\begin{equation} \label{eq:329}
    V^{[1]}(t) = V^{[2]}(t) \quad \forall t\in [0, T].
\end{equation}

\subsubsection{Uniqueness of \texorpdfstring{$m$}{m} in the general case}
To see that $m^{[1]} = m^{[2]}$, observe that the second term in \eqref{eq:323} is identically zero.
\begin{equation} \label{eq:339a}
    \dfrac{1}{2(2 + \epsilon \overline{\eta})} \Big [ \Big (\sum_{k=1}^{K} \delta m_k (1 - \overline{\Delta V_k})\Big) + \epsilon \overline{\phi_\epsilon} \delta \eta\Big]^2 = 0 .
\end{equation}
With uniqueness of $V$ in mind, combining \eqref{eq:339a} with \eqref{eq:321a} results in $\delta\phi_\epsilon = 0$. Therefore, $\phi_{\epsilon}^{[1]}(t) = \phi_{\epsilon}^{[2]}(t)$  for all $t\in [0, T]$. Consequently, the equations in \eqref{eq:317a} become
\begin{equation} \label{eq:333}
\begin{split}
    &\dfrac{d \delta m_k}{dt} = -  G_k \delta m_k +G_{k+1} \delta m_{k+1}, \quad k \in \mathcal{K}_{-1}\\
    &\dfrac{d \delta m_K}{dt} = - G_{K} \delta m_K.
\end{split}
\end{equation}
Since $\delta m_K (0) = 0$, the solution when $k = K$ is $\delta m_K (t) = 0$ for all $t\in [0, T]$. By backward induction, we find that for any $ k \in \mathcal{K}$, we have $\delta m_k (t) = 0 $ for all $t\in[0, T]$. Therefore, 
\begin{equation} \label{eq:336}
    m^{[1]}(t) = m^{[2]}(t) \quad \forall t\in [0, T].
\end{equation}
This concludes the proof of uniqueness. \qed

\bibliographystyle{plainnat} 
\bibliography{references}

\end{document}